\documentclass[a4paper,11pt]{article}

\usepackage{xcolor}
\usepackage[pdfencoding=unicode, hidelinks]{hyperref}
\hypersetup{
	colorlinks,
	linkcolor={red},
	citecolor={blue},
}

\usepackage[english]{babel}
\usepackage[utf8]{inputenc}
\usepackage{amsthm}
\usepackage{amsmath}
\usepackage{amssymb}
\usepackage{mathrsfs}
\usepackage[mathscr]{eucal}
\usepackage{wasysym}
\usepackage{textcomp}
\usepackage{enumerate}
\usepackage{setspace}
\usepackage{latexsym}
\usepackage{graphicx}
\usepackage{mathtools}
\usepackage{tabularray}
\usepackage{geometry}
\usepackage{float}
\usepackage{emptypage}
\usepackage{commath}
\usepackage{multicol}
\usepackage{titlesec}
\usepackage[shortlabels]{enumitem}
\usepackage{verbatim}
\usepackage{cases}
\usepackage{cleveref}
\usepackage{orcidlink}
\usepackage{authblk}
\usepackage{mathbbol}

\usepackage{csquotes}
\usepackage[backend=biber,
    style=alphabetic,
    citestyle=alphabetic,
    sorting=nyt, 
    doi=false,
    url=false,
    eprint=false,
    isbn=false,
    maxnames=6,
    abbreviate=true,
    giveninits=true]{biblatex}
\usepackage{hyperref}
\AtEveryBibitem{%
  \clearfield{note}%
}

\theoremstyle{plain}
\newtheorem{thm}{Theorem}[section]
\theoremstyle{plain}

\theoremstyle{plain}

\theoremstyle{plain}
\newtheorem{defn}[thm]{Definition}
\theoremstyle{plain}

\theoremstyle{definition}
\newtheorem{remark}[thm]{Remark}

\numberwithin{equation}{section}

\newcommand{\RR}{\mathbb{R}}
\newcommand{\NN}{\mathbb{N}}
\newcommand{\dx}{\,\ensuremath{\mathrm{d}}x}
\newcommand{\ds}{\,\ensuremath{\mathrm{d}}s}
\newcommand{\dt}{\,\ensuremath{\mathrm{d}}t}

\renewcommand{\norm}[1]{\lVert#1\rVert}
\renewcommand{\phi}{\varphi}
\newcommand{\duality}[2]{\langle#1,#2\rangle}
\newcommand{\dualprod}[2]{\langle#1,#2\rangle_{\Vs}}
\newcommand{\intprod}[2]{(#1,#2)_{\Hs}}

\newcommand{\Ws}{W}
 
\newcommand{\Vs}{V}
\newcommand{\Hs}{H}

\newcommand{\nnu}{\boldsymbol{\nu}}

\newcommand{\R}{\mathcal{R}}
\newcommand{\N}{\mathcal{N}}

\newcommand{\hbeta}{\widehat{\beta}}
\newcommand{\hbetae}{\widehat{\beta}_{\varepsilon}}
\newcommand{\betae}{\beta_{\varepsilon}}
\newcommand{\hpi}{\widehat{\pi}}
\newcommand{\hf}{\widehat{f}}

\newcommand{\lambdap}{\lambda_P}
\newcommand{\lambdaa}{\lambda_A}
\newcommand{\lambdae}{\lambda_E}
\newcommand{\lambdac}{\lambda_C}
\newcommand{\lambdab}{\lambda_B}
\newcommand{\lambdad}{\lambda_D}

\DeclareSymbolFontAlphabet{\mathbb}{AMSb}
\DeclareSymbolFontAlphabet{\mathbbm}{bbold}
\newcommand{\h}{\mathbbm{h}}
\renewcommand{\k}{\mathbbm{k}}

\newcommand{\difft}{\frac{\text{d}}{\text{d}t}}
\DeclareMathOperator{\diver}{div}
\newcommand{\into}{\int_{\Omega}}
\newcommand{\intTo}{\int_0^T \!\!\!\!\into}
\newcommand{\intTno}{\int_0^{T^n}\!\!\!\!\!\into}
\newcommand{\intTn}{\int_0^{T^n}}
\newcommand{\intto}{\int_0^t \!\!\!\into}
\newcommand{\intt}{\int_0^t}

\newcommand{\thetan}{\theta^n}
\newcommand{\thetanj}{\theta^n_j}
\newcommand{\thetae}{\theta_{\varepsilon}}
\newcommand{\phin}{\phi^n}
\newcommand{\phinj}{\phi^n_j}
\newcommand{\ophin}{\overline{\phi}^n}
\newcommand{\phie}{\phi_{\varepsilon}}

\newcommand{\mun}{\mu^n}
\newcommand{\munj}{\mu^n_j}
\newcommand{\omun}{\overline{\mu}^n}
\newcommand{\mue}{\mu_{\varepsilon}}
\newcommand{\sigman}{\sigma^n}
\newcommand{\sigmanj}{\sigma^n_j}
\newcommand{\sigmae}{\sigma_{\varepsilon}}

\crefname{lemma}{lemma}{lemmas}
\Crefname{lemma}{Lemma}{Lemmas}
\crefname{prop}{proposition}{proposition}
\Crefname{prop}{Proposition}{Propositions}
\crefname{cor}{corollary}{corollaries}
\Crefname{cor}{Corollary}{Corollaries}
\crefname{remark}{remark}{remarks}
\Crefname{remark}{Remark}{Remarks}
\crefname{thm}{theorem}{theorems}
\Crefname{thm}{Theorem}{Theorems}
\Crefname{section}{Section}{Sections}

\begin{document}

\begin{center}
		
    {\Large \bf Well-posedness for a fourth-order nonisothermal\\[2mm]
    tumor growth model of Caginalp type}

    \vskip0.5cm
    
    {\large\textsc{Giulia Cavalleri\orcidlink{0009-0006-4154-9659}$^1$}} \\
    {\normalsize e-mail: \texttt{giulia.cavalleri01@universitadipavia.it}} \\
    \vskip0.35cm

    {\large\textsc{Pierluigi Colli\orcidlink{0000-0002-7921-5041}$^{1,2}$}} \\
    {\normalsize e-mail: \texttt{pierluigi.colli@unipv.it}} \\
    \vskip0.35cm

    {\large\textsc{Elisabetta Rocca\orcidlink{/0000-0002-9930-907X}$^{1,2}$}} \\
    {\normalsize e-mail: \texttt{elisabetta.rocca@unipv.it}} \\
    \vskip0.35cm
    
    {\footnotesize $^1$Department of Mathematics ``F. Casorati'', University of Pavia, 27100 Pavia, Italy}
    \vskip0.1cm
    
    {\footnotesize$^2$Research Associate at the IMATI -- C.N.R. Pavia, via Ferrata 5, 27100 Pavia, Italy}
    \vskip0.5cm
		
\end{center}

\begin{abstract}
    \noindent We introduce a nonisothermal phase-field system of Caginalp type that describes tumor growth under hyperthermia. The model couples a possibly viscous Cahn--Hilliard equation, governing the evolution of the healthy and tumor phases, with an equation for the heat balance, and a reaction-diffusion equation for the nutrient concentration.  The resulting nonlinear system incorporates chemotaxis and active transport effects, and is supplemented with no-flux boundary conditions. The analysis is carried out through a two-step approximation procedure, involving a regularization of the potential and a Faedo--Galerkin discretization scheme. Under stronger regularity assumptions, we further establish the existence of strong solutions and their uniqueness via a continuous dependence result.

    \vskip3mm
    
    \noindent {\bf Key words:} nonlinear initial-boundary value problem, Cahn--Hilliard equation, nonisothermal model, well-posedness, regularity results, tumor growth.

    \vskip3mm
    
    \noindent {\bf AMS (MOS) Subject Classification:} 
    35K61, 
    35K35, 
    35D30, 
    35B65, 
    35A02, 
    35Q92, 
    92C50. 

\end{abstract}

\section{Introduction}
Despite the advances in modern medicine, cancer remains one of the leading causes of death worldwide (see, e.g., \cite{Ferlay_etal}). 
Moreover, due to the actual rates in population growth and aging, the number of new cases per year is set to increase, further intensifying the global burden of the disease (see, e.g., \cite{Bray_etal}). As a result, there is a great interest in combining existing therapeutic strategies to enhance treatment efficacy and improve patient outcomes. Mathematical modeling plays a crucial role in this context, providing a powerful framework for exploring the complex interplay between biological processes and different therapies. 
In particular, it enables a systematic and quantitative investigation of treatment combinations, dosages, and scheduling, which is not always possible in clinical trials, providing valuable insights. 
In this work, we focus on the effects of hyperthermia therapy, a treatment approach that has been explored for several decades (see, e.g., \cite{Mallory_etal} and the references therein). Hyperthermia consists of raising the local, regional, or whole-body temperature above 39°C for 30 to 60 minutes, using techniques such as radio wave, microwave, or focused ultrasound heating, perfusion, and thermal chambers (see, e.g., \cite{Habash_etal_2011}).
It is sometimes used alone, but most frequently is combined with other primary therapies, such as chemotherapy, radiotherapy, and immunotherapy. Depending on the temperature, it has different effects.
\begin{enumerate}[(i)]
    \item Moderate hyperthermia (below 42°C) increases tumor perfusion, primarily due to heat-induced vasodilatation. This, in turn, improves the delivery of chemotherapeutic or immunotherapeutic drugs.
    \item High hyperthermia (i.e., temperatures between 42°C and 50°C) exerts direct cytotoxic effects and induces vascular damage. Moreover, it impairs DNA repair mechanisms, thereby potentially enhancing the susceptibility of tumor cells to other treatments such as chemotherapy or radiotherapy.
    \item Thermal ablation (above 50°C), applied directly in situ, results in irreparable cellular damage and consequent apoptosis and necrosis of tumor tissue. 
\end{enumerate}
The result depends on several factors, including treatment timing, heat, and nutrient distribution within the tissue. To address these questions, we derive a model that couples the dynamics of tumor growth, nutrient 
diffusion, and localized thermal treatment. We aim to rigorously investigate the analytical properties of the resulting system, such as the existence of weak solutions, regularity, and continuous dependence on the problem's data.\\

\noindent \textbf{The PDE system.} We consider a smooth domain $\Omega \subseteq \RR^d$ with $d = 2,3$ and final time $T>0$. The PDE system reads as
\allowdisplaybreaks
\begin{subequations}\label{eq:problem}
    \begin{align}
        & \partial_t(\theta +\ell\phi) - \Delta \theta = u,\label{eq:temperature}\\
        & \partial_t\phi - \Delta \mu = (\lambdap \sigma - \lambdaa - \lambdae \theta)\h(\phi),\label{eq:ch1}\\
        & \mu = \tau \partial_t \phi - \Delta \phi + \beta(\phi) + \pi(\phi) - \chi \sigma -\Lambda \theta, \label{eq:ch2}\\
        & \partial_t \sigma - \Delta (\sigma - \chi \phi) = - \lambdac \sigma \h(\phi) + \lambdab (\sigma_B - \sigma) - \lambdad \sigma \k(\theta), \label{eq:nutrient}
    \end{align}
\end{subequations}
in the parabolic cylinder $Q \coloneqq \Omega \times (0,T)$, coupled with the boundary conditions
\begin{equation}\label{eq:boundary_cond}
    \partial_{\nnu} \theta = \partial_{\nnu} \phi = \partial_{\nnu} \mu = \partial_{\nnu} \sigma = 0
\end{equation}
on $\Sigma \coloneqq \partial \Omega \times (0,T)$, and with the initial conditions
\begin{equation}\label{eq:initial_datum}
     \theta(0) = \theta_0, \quad \phi(0) = \phi_0, \quad \sigma(0) = \sigma_0
\end{equation}
in $\Omega$. 
The subsystem \eqref{eq:temperature}--\eqref{eq:ch2} is of Caginalp type, and is made of the second-order parabolic equation \eqref{eq:temperature}, ruling the evolution of $\theta$---the relative temperature with respect to a certain critical value normalized to 0---and the (possibly) viscous Cahn--Hilliard equation \eqref{eq:ch1}--\eqref{eq:ch2}, governing the evolution of the phase-field variable $\phi$, which is the difference in volume fraction between tumor and healthy cells.

Within this framework, in place of a sharp interface separating tumor and healthy tissues, there is a narrow transition layer, where both cell types coexist. The result is a model well-suited for describing topological changes such as coalescence or breaking up phenomena, which frequently occur during the early stages of tumor development. 
From the mathematical point of view, at least in principle, $\phi$ has values in the interval $[-1,1]$, where $\{\phi = 1\}$ identifies the tumor region, $\{\phi = -1\}$ the healthy region, and $\{-1 < \phi < 1\}$ the diffuse interface (see \Cref{remark:double_well_potential}).

Finally, the parabolic equation \eqref{eq:nutrient} prescribes the dynamics of the nutrient concentration $\sigma$, a chemical species dissolved in the blood (such as oxygen or glucose) and consumed only by tumor cells. 

Before delving into the derivation of the model introduced above, let us recall some of the related literature. The widely studied and now well-understood Cahn--Hilliard equation is usually set when assuming constant temperature, which is not always reasonable from the modeling point of view (see \cite{Miranville_19} for further details). Therefore, a number of nonisothermal phase transition models have been proposed (see, e.g., \cite{Caginalp_88, Penrose_Fife_90, Alt_Pawlow_92} and the more recent \cite{Miranville_Schimperna_05, Marveggio_Schimperna_21, DeAnna_etal_24}). 
In particular, the Caginalp model has been discussed and analyzed over the years, also recently: on that, we may quote~\cite{Colli_etal_23, Colli_etal_24} and references therein.

As for tumor growth models based on the Cahn--Hilliard equation, the literature is extensive  (see, e.g., \cite{Frigeri_Grasselli_Rocca_15, Garcke_Lam_Sitka_Styles_16, Garcke_Lam_17, Ebenbeck_Garcke_Nurnberg_21}). The nutrient dynamic is typically taken into account, but various works also incorporate the evolution of other biological quantities that may affect cancer 
development, such as relevant biomarkers,  elastic effects, and fluid flows (see, e.g., \cite{Colli_Gomez_Lorenzo_etal_20, Cavalleri_Colli_Miranville_Rocca_2026, Garcke_Lam_Signori_21, Garcke_Kovacs_Trautwein_22, Garcke_Lam_Sitka_Styles_16, Knopf_Signori_22} and the references therein). 

However, to the best of our knowledge, only the recent contribution \cite{Ipocoana_22} incorporates temperature effects. In the paper~\cite{Ipocoana_22}, a nonisothermal diffuse interface model for tumor growth is proposed and analyzed. The system tries to capture the complex coupling between temperature variation, nutrient transport, cellular proliferation, and apoptosis within a Cahn--Hilliard-type framework. The model is formulated in terms of an entropy balance and is rigorously studied in a suitable weak setting. In particular, the author establishes the existence of solutions for the associated initial-boundary value problem. 
While our approach shares with~\cite{Ipocoana_22} the interest in coupling phase-field dynamics with thermal effects and a tumor growth process, we address a different system and adopt a distinct analytical strategy. 
In particular, beyond proving the existence of weak solutions, we also establish higher regularity and continuous dependence on the data. These analytical results lay the groundwork for the analysis of an associated optimal control problem, to be addressed in a forthcoming work~\cite{Cavalleri_Colli_Rocca_inprep}.\\

\noindent \textbf{Derivation of the model.} The model~\eqref{eq:problem}--\eqref{eq:initial_datum} is derived following the approach from Caginalp (see 
\cite{Caginalp_88,Caginalp_90} but also \cite[][Chapter 4, Example 4.4.2]{Brokate_Sprekels_96}). The evolution of temperature is governed by the heat balance equation
\begin{equation}
    \partial_t (\theta + \ell \phi) -\Delta \theta = u,
\end{equation}
where $u$ is a given heat source, and $\ell$ is a constant related to the tissue's latent heat. This equation is coupled with the mass balance 
\begin{equation}\label{eq:phi_balance_law}
    \partial_t\phi + \diver(- \nabla \mu) = U^{\phi},
\end{equation}
where $U^{\phi}$ denotes a phenomenologically chosen source term, and $\mu$ is the chemical potential of the system. Following the formulation in~\cite{NovickCohen_88}, the chemical potential is defined as
\begin{equation}
    \mu \coloneqq \tau \partial_t \phi + \frac{\delta \mathcal{F}}{\delta \phi}(\theta,\phi,\sigma), 
\end{equation}
where the first addend provides a viscous, dissipative contribution with $\tau \geq 0$, while the second term is the variational derivative of the total free energy functional, which is assumed to be given by
\begin{equation}\label{eq:energy}
    \mathcal{F}(\theta,\phi,\sigma) = \into   \left[\frac{1}{2}|\nabla \phi|^2 +  \hbeta(\phi) + \hpi(\phi) + \frac{1}{2}|\sigma|^2 + \chi \sigma (1 -\phi) -\Lambda \theta \phi \right]\dx.
\end{equation}
The term $\frac{1}{2} |\nabla \phi|^2 + \hbeta(\phi) + \hpi(\phi)$ is of classical Ginzburg--Landau type: the gradient accounts for the interfacial energy, while the remaining part is a double-well potential related to the interaction between the two phases. 
The term $\frac{1}{2} |\sigma|^2$ reflects the energetic cost associated with the presence of nutrients, implying that high concentrations increase the system’s free energy.  
The two final addends couple the Cahn--Hilliard equation with the nutrient and with the temperature equations. The nonnegative constants $\chi$ and $\Lambda$ represent, respectively, a transport coefficient (modeling effects such as chemotaxis and active transport) and a parameter related to the tissue's latent heat. The nutrient equation is also obtained by a mass balance, choosing the flux as the variational derivative of the free energy $\mathcal{F}$ with respect to $\sigma$:
\begin{equation}\label{eq:sigma_balance_law}
    \partial_t \sigma + \diver \bigg[-\nabla \bigg(\frac{\delta \mathcal{F}}{\delta \sigma}(\theta,\phi,\sigma)\bigg) \bigg] = U^{\sigma},
\end{equation}
for a nutrient source $U^{\sigma}$, again phenomenologically chosen. 
We observe that the Cahn--Hilliard equation \eqref{eq:phi_balance_law} includes the term $\diver (\chi \nabla \sigma)$, accounting for chemotaxis, i.e., the tendency of tumor cells to migrate toward regions with higher nutrient concentrations. Additionally, the nutrient equation \eqref{eq:sigma_balance_law} contains the term $-\diver(-\chi \nabla \phi)$, accounting for the so-called active transport, which models the preferential flow of nutrients toward tumor cells. This is motivated by the biological observation that tumor cells often overexpress specific membrane transporters facilitating nutrient uptake (for more details, see the \cite{Garcke_Lam_Sitka_Styles_16,Ebenbeck_Garcke_Nurnberg_21}). 

Finally, we introduce the mass and nutrient sources that we are going to employ throughout the paper. Starting from $U^{\phi}$, we assume 
\begin{equation*}
    U^{\phi}(\theta,\phi,\sigma) \coloneqq (\lambdap \sigma - \lambdaa - \lambdae \theta)\h(\phi).
\end{equation*}
According to the literature (see \cite{Garcke_Lam_17}), we assume the mechanisms controlling cell division to be suppressed in tumor cells. Thus, proliferation is limited only by the availability of nutrients. We model it with the term $\lambda_p \sigma$, where $\lambda_p$ is a fixed proliferation coefficient. We also suppose that tumor cells only die because of apoptosis, and we denote with $\lambda_a$ the constant apoptosis rate. Moreover, here we assume that $\theta$ has cytotoxic effects proportional to the temperature, which we incorporate through the term $-\lambdae \theta$, where $\lambdae$ is a fixed positive parameter. The function $\h$ guarantees that the phenomena we have just described are proportional to the tumor cells available in a certain area. For example, $\h$ may be a monotone increasing function which is 0 where $\{\phi = -1\}$ and 1 where $\{\phi = 1\}$. On the other hand, we assume the nutrient source $U^{\sigma}$ to be of the form
\begin{equation*}
     U^{\sigma}(\theta,\phi,\sigma) \coloneqq - \lambdac \sigma \h(\phi) + \lambdab (\sigma_B - \sigma) - \lambdad \sigma \k(\theta).
\end{equation*}
The first two addends are again compliant with \cite{Garcke_Lam_17}. The term $- \lambdac \sigma \h(\phi)$ models the fact that tumor consumption is higher where tumor cell density is higher, while $\lambdab (\sigma_B - \sigma)$ accounts for nutrient supply from the preexisting capillaries. A novel aspect of the model is the inclusion of the term~$-\lambdad \sigma \k(\theta)$, which describes the influence of temperature on nutrient absorption---e.g., due to vasodilation effects that enhance nutrient transport in warmer tissue regions (see, e.g., \cite{Song_etal_05}).
We couple the system~\eqref{eq:problem} hereby obtained with the no-flux boundary conditions~\eqref{eq:boundary_cond}, which  reflect the idealization that the system is effectively isolated from its surroundings, and with the initial 
conditions~\eqref{eq:initial_datum}.\\

\noindent \textbf{Plan of the paper.} The paper is organized as follows. In \Cref{sect:notation_preliminaries}, we introduce the notation and some mathematical preliminaries we use throughout this work. In \Cref{sect:hyp_main_result}, we state our main theorems and the hypotheses we assume. Then, we prove the previously mentioned results. More precisely, \Cref{sect:existence} is devoted to the existence of weak solutions of the system \eqref{eq:problem}--\eqref{eq:initial_datum}. \Cref{sect:reg} regards the additional regularity we are able to reach, strengthening the assumptions on the initial data. Finally, in \Cref{sect:cont_dep} we prove a continuous dependence result which, in turn, implies the uniqueness of strong solutions. 

\section{Notation and preliminaries}\label{sect:notation_preliminaries}
\textbf{Notation.} For any Banach space $(X,\norm{\cdot}_{X})$ we employ $(X^*,\norm{\cdot}_{X^*})$ for its topological dual. We denote the Lebesgue and Sobolev spaces over $\Omega$ respectively as $L^p(\Omega)$ and $W^{k,p}(\Omega)$. In the special case of $p=2$, we set $H^k(\Omega) \coloneqq W^{k,2}(\Omega)$. The norm of the Bochner space $W^{k,p}(0,T;X)$ is indicated as $\norm{\cdot}_{W^{k,p}(X)}$, omitting the time interval $(0,T)$ for the sake of brevity. If the final time differs from $T$, it will be written explicitly to avoid ambiguity. With the notation $C^0([0,T];X)$, we mean the space of continuous $X$-valued functions, while $C^0_w([0,T];X)$ denotes the space of functions $v\in L^\infty (0,T;X)$ that are weakly continuous from $[0,T]$ to $X$.
For convenience, we introduce the notation
\begin{equation*}
    W \coloneqq \{v \in H^2(\Omega) \, | \, \partial_{\nnu} v = 0 \text{ on } \partial \Omega\}, \qquad \Vs \coloneqq H^1(\Omega), \qquad \Hs \coloneqq L^2(\Omega).
\end{equation*}
We use $\intprod{\cdot}{\cdot}$ for the internal product in $\Hs$ and $\dualprod{\cdot}{\cdot}$ for the dual pairing between $\Vs^*$ and $\Vs$. As always in this setting, we identify $\Hs$ with its dual space $\Hs^*$ and as a subspace of $V^*$ through
\begin{equation*}
    \dualprod{w}{v} = \intprod{w}{v}
\end{equation*}
for every $w \in \Hs$ and  $v \in \Vs$. This way we obtain the usual Hilbert triplet 
\begin{equation*}
     \Vs \subset \subset \Hs \subset \subset  \Vs^*,
\end{equation*}
where the embeddings are dense and compact. For every $v^* \in \Vs^*$, we introduce its (generalized) mean value as 
\begin{equation*}
    \overline{v^*}\coloneqq \frac{1}{|\Omega|} \dualprod{v^*}{1}
\end{equation*}
for every $v^* \in \Vs^*$. Obviously, it coincides with the usual mean value if $v \in \Hs$ and, by extension, if $v \in L^1(\Omega)$. We employ
\begin{equation*}
   \Vs_0, \qquad \Hs_0, \qquad \Vs_0^*
\end{equation*}
for the closed subspace respectively of $\Vs$, $\Hs$, $\Vs^*$ of elements with zero mean value.\\ 

\noindent \textbf{The Neumann--Laplace operator.} We introduce the $- \Delta$ operator with Neumann boundary condition restricted to the space $\Vs_0$ as 
\begin{align*}
    &\R \,: \, \Vs_0 \to \Vs_0^* \text{ s.t.}\\
    &\dualprod{\R v }{w} \coloneqq \into \nabla v \cdot \nabla w \dx \quad \forall w \in \Vs_0.
\end{align*}
$\R$ is an isomorphism, and we denote its inverse by
\begin{align*}
    \N \coloneqq \R^{-1} \,: \, \Vs_0^* \to \Vs_0.
\end{align*}
The operators $\R$ and $\N$ enjoy several well-known properties, which we enlist here for the reader's convenience: 
\begin{align}
    & \text{$\dualprod{v^*}{\N w^*} = \dualprod{w^*}{\N v^*} =  \intprod{\nabla \N w^*}{\nabla \N v^*}$\quad for all $v^*, w^* \in \Vs_0^*$,}\\
    & \text{$ \intprod{w}{v}= \intprod{w-\overline{w}}{v} =  \into \nabla w \cdot \nabla \N v \dx$\quad for all $w \in \Vs, \ v \in \Hs_0$,} \label{enum:eqN}\\
    & \text{$\norm{v}_{\Hs} \leq  \norm{\nabla v}_{\Hs}^{1/2} \norm{\nabla \N v}_{\Hs}^{1/2} $\quad for every $v \in \Vs_0$ \label{enum:int_ineq_N},}\\
    \begin{split}
         &\text{$\displaystyle \dualprod{\partial_t v^*(t)}{\N v^*(t)} = \frac{1}{2} \difft \norm{\nabla \N v^*(t)}_{\Hs}^2$}\\ 
         &\quad\quad\text{for a.e. $t \in (0,T)$, for all  $v^* \in H^1(0,T; \Vs^*_0)$.}
    \end{split}
    \end{align}
Since $\N$ is defined over the space $\Vs_0^*$ but the functions we work with do not have, in general, zero mean value, it is useful to notice that the standard norms over $\Ws$, $\Vs$, $\Hs$, and $\Vs^*$ are equivalent to the following ones:
\begin{alignat*}{2}
     &\norm{v}_{\Ws}  &&\simeq (\norm{-\Delta v}_{\Hs}^2 +  \overline{v}^2)^{1/2},\\
     &\norm{v}_{\Vs}  &&\simeq (\norm{\nabla v}_{\Hs}^2 + \overline{v}^2)^{1/2},\\
     &\norm{v}_{\Hs}  &&\simeq (\norm{v - \overline{v}}_{\Hs}^2 + \overline{v}^2)^{1/2},\\
     &\norm{v}_{\Vs^*} &&\simeq (\norm{v -  \overline{v}}_{\Vs^*}^2 + \overline{v}^2)^{1/2} \simeq (\norm{\nabla \N (v- \overline{v})}_{\Hs}^2 + \overline{v}^2)^{1/2}.
\end{alignat*}

\addvspace{0.6cm}
\noindent \textbf{Useful inequalities.} In our proofs, we frequently rely on a set of classical inequalities. In particular, we make repeated use of the Poincaré, the H\"older, and the Young inequalities, the latter of which reads
\begin{equation*}
    ab \leq \eta a^2 + \frac{1}{4\eta}b^2 \quad \text{for every } a,b \in \RR,\ \eta > 0.
\end{equation*}
Another important tool is Ehrling's Lemma (see, e.g., \cite[][Theorem 16.4, p. 102]{Lions_Magenes_12}), a well-known compactness result which states that, given three Banach spaces $X \subset \subset Y \subset Z$, for every $\eta > 0$ there exists a $C_{\eta} > 0$ such that
\begin{equation*}
    \norm{v}_Y \leq \eta \norm{v}_X + C_{\eta} \norm{v}_Z
\end{equation*}
for every $v \in X$. In our setting, we apply it for $\Vs \subset\subset L^p(\Omega) \subset \Vs^*$ and $p=2,4$, which is true in dimension $d=2,3$.\\

\noindent Finally, $C$ is a constant depending only on the problem's data and whose value might change from line to line or even within the same line. If we want to highlight a dependency on a certain parameter, we put it as a subscript (e.g., $C_{\tau}$ indicates a constant that depends on $\tau$, $C_0$ a constant that depends on the initial data, etc.).

\section{Hypotheses and main results}\label{sect:hyp_main_result}
Throughout the paper, we consider a bounded domain $\Omega$  of class $C^2$ in $\RR^d$ with $d = 2,3$, and we fix a final time $T > 0$.

\begin{enumerate}[\rm{(H\arabic*)}]
\item \label{hyp:constants}
We assume that 
\begin{align}
    & \ell, \Lambda, \chi \text{ are positive constants},\\
    & \tau \text{ and }\lambdap, \lambdaa, \lambdae, \lambdac, \lambdab, \lambdad  \text{ are real nonnegative constants}.
\end{align}
\item \label{hyp:given_function}
The assigned functions $u$ and $\sigma_B$ enjoy the regularity
\begin{align}
    & u \in L^{\infty}(Q) \ \text{ with } \ \norm{u}_{L^{\infty}} \leq M, \label{hyp:control}\\
    & \sigma_B \in L^2(Q)
\end{align}
for a nonnegative constant $M$. 

\item \label{hyp:nonlinearities}
Regarding the nonlinearities $\h$ and $\k$, we require that
\begin{align}
    & \h, \k \in C^{0,1}(\RR), \text{ and}\\
    &0 \leq \h \leq \h^*,\ 0 \leq \k \leq \k^*,
\end{align}
for two nonnegative constants $\h^*,\ \k^*$.
\item \label{hyp:beta_pi}
We consider a potential $\hf = \hbeta + \hpi$ split into the sum of a convex part and a nonconvex perturbation. Explicitly, we suppose
\begin{align}
     & \hbeta, \,  \hpi \in C^2(\RR), \label{pier19} \\
     & \hbeta \text{ is convex and nonnegative, with }
     \hbeta (0)=0. \label{pier20}
\end{align}
Moreover, denoting its derivative as $f \coloneqq \beta + \pi$ where $\beta \coloneqq \hbeta'$ and $\pi \coloneqq \hpi'$, the following growth conditions hold:
\begin{align}
    & |\beta(r)| \leq C_{\beta} (\hbeta(r) + 1),
    \label{hyp:growth_beta}\\[1mm]
    & |\pi'(r)| \leq C_{\pi}
    \label{hyp:pi_lip}
\end{align}
for all $r \in \RR$, where $C_{\beta}$, $C_{\pi}$ are given nonnegative constants. We also point out that
\begin{equation}
   \beta(0) = 0 .
    \label{pier-beta}
\end{equation}
\end{enumerate}

\begin{remark}\label{remark:double_well_potential}
    It follows that $\hpi$ has at most quadratic growth and $\hbeta$ has at most exponential growth. A meaningful example covered by our hypotheses is the regular quartic potential
    \begin{equation*}
        \hf_{reg}(r) = \frac{1}{4}(r^2 - 1)^2.
    \end{equation*}
    Notice that the logarithmic potential, as well as the double obstacle potential, is instead excluded. In particular, we are not able to guarantee that the order parameter $\phi$ assumes values in the physically relevant interval $[-1,1]$.
\end{remark}

\begin{remark}\label{remark:pier}
We observe that the assumption that $\widehat{\beta}$ attains its minimum value $0$ at $0$ (see~\eqref{pier20} and~\eqref{pier-beta}) is not restrictive. Indeed, if this were not the case, one could modify $\widehat{\beta}$ by subtracting its tangent line at $0$, so that the new function $\widehat{\beta}$ remains convex and nonnegative, and additionally attains its minimum $0$ at $0$.
Moreover, there is no issue in adding a linear term to $\widehat{\pi}$; doing so still preserves the Lipschitz continuity of $\pi$.
\end{remark}

\begin{defn}\label{defn:weak_sol}
    We define a weak solution of the PDE system~\eqref{eq:problem}--\eqref{eq:initial_datum} as a quadruplet $(\theta, \phi, \mu, \sigma)$ with the regularity
    \begin{gather}
        \theta \in  H^1(0,T;\Vs^*) \cap  C^0([0,T];\Hs) \cap  L^2(0,T;\Vs),\\  
        \phi \in   H^1(0,T;\Vs^*) \cap  C^0_w([0,T];\Vs) \cap  L^2(0,T;\Ws) ,\\
        \mu \in L^2(0,T; \Vs),\\
        \sigma \in  H^1(0,T;\Hs) \cap  C^0([0,T];\Vs) \cap  L^2(0,T;\Ws), 
    \end{gather}
    such that equations \eqref{eq:temperature}--\eqref{eq:nutrient} and the boundary conditions \eqref{eq:boundary_cond} are satisfied in the following variational sense
    {
    \allowdisplaybreaks
    \begin{subequations}\label{eq:problem_eq_with_spaces}
        \begin{align}
            &\duality{\partial_t(\theta +\ell
        \phi)}{v}_{\Vs} + \into  \nabla \theta \cdot \nabla v \dx = \into u v \dx, \label{eq:temperature_def}\\
            &\duality{\partial_t \phi}{v}_{\Vs} + \into  \nabla \mu \cdot \nabla v \dx = \into (\lambdap \sigma - \lambdaa - \lambdae \theta) \h(\phi) v \dx, \label{eq:ch1_def}\\
            & \begin{aligned}
            	&\duality{\tau \partial_t \phi}{v}_{\Vs} + \into \nabla \phi\cdot \nabla v \dx\\
            	& \quad \quad + \into \left(\beta(\phi) + \pi(\phi) - \chi \sigma - \Lambda\theta\right) v \dx  = \into \mu v \dx, 
            \end{aligned} \label{eq:ch2_def} \\
            & \begin{aligned}
            	&\duality{\partial_t\sigma}{v}_{\Vs} + \into \nabla \sigma \cdot \nabla v \dx  - \chi \into \nabla \phi \cdot \nabla v \dx\\
            	&\quad \quad = \into \left(-\lambdac \sigma \h(\phi) + \lambdab (\sigma_B - \sigma) - \lambdad \sigma \k(\theta)\right) v \dx,
            \end{aligned} \label{eq:nutrient_def}
        \end{align}
    \end{subequations}
	}
    a.e. in $(0,T)$ for all $v \in \Vs$, and the initial data \eqref{eq:initial_datum}
    \begin{equation}
         \theta(0) = \theta_0, \quad \phi(0) = \phi_0, \quad  \sigma(0) = \sigma_0
    \end{equation}
    a.e. in $\Omega$.
\end{defn}

\begin{remark}
An equivalent definition of weak solution is obtained by integrating in time the equalities in~\eqref{eq:problem_eq_with_spaces} over the interval $(0,T)$ and choosing test functions $v \in L^2(0,T; \Vs)$. In addition, we point out that \eqref{eq:ch2_def} and \eqref{eq:nutrient_def} can be equivalently rewritten as equations that hold almost everywhere and are complemented by homogeneous Neumann boundary conditions as in~\eqref{eq:boundary_cond}.
\end{remark}

\begin{thm}[Existence of weak solutions]\label{thm:existence}
    Assume that the set of hypotheses \ref{hyp:constants}--\ref{hyp:beta_pi} holds and that the initial data satisfy
    \begin{equation}\label{hyp:initial_condition}
        \theta_0 \in \Hs,\quad \phi_0 \in \Vs,\quad \hbeta(\phi_0) \in L^1(\Omega),\quad \sigma_0 \in \Vs.
    \end{equation}
    Then, the PDE system~\eqref{eq:problem}--\eqref{eq:initial_datum} admits at least a weak solution $(\theta, \phi, \mu, \sigma)$ with the additional regularity
    \begin{equation}
        \tau^{1/2} \phi \in H^1(0,T; \Hs),
    \end{equation}
   which satisfies the estimate
    \begin{equation}
        \begin{split}
            &\norm{\theta}_{H^1(\Vs^*) \cap L^{\infty}(\Hs) \cap L^2(\Vs)} + \norm{\phi}_{H^1(\Vs^*) \cap L^{\infty}(\Vs) \cap L^2(\Ws)} + \norm{\tau^{1/2}\phi}_{H^1(\Hs)}\\
            &\quad + \norm{\beta(\phi)}_{L^2(\Hs)} + \norm{\mu}_{L^2(\Vs)} + \norm{\sigma}_{H^1(\Hs) \cap L^{\infty}(\Vs) \cap L^2(\Ws)}\leq C_1
        \end{split}
        \label{est-pier1}
    \end{equation}
    for a constant $C_1>0$ that depends on $M$ and on the other problem's data. 
\end{thm}

\noindent If we improve the regularity of the initial data, we prove that the system has a more regular solution, according to the following result.

\begin{thm}[Regularity]\label{thm:regularity}
    Assume that the hypotheses \ref{hyp:constants}--\ref{hyp:beta_pi} are satisfied, and that the initial data enjoy
    \begin{equation}\label{hyp:reg_initial_condition}
        \theta_0 \in \Vs \cap L^{\infty}(\Omega),\quad \phi_0 \in \Ws \cap H^3(\Omega),
        \quad \sigma_0 \in \Vs.
    \end{equation} 
    Then, there exists a weak solution $(\theta, \phi, \mu, \sigma)$ with the additional regularity
        \begin{gather}
        \theta \in  H^1(0,T;\Hs) \cap L^\infty(0,T;\Vs)  \cap L^2(0,T;\Ws) \cap L^{\infty}(Q),\\  
        \phi \in   W^{1,\infty}(0,T;\Vs^*) \cap  H^1(0,T;\Vs) \cap  L^{\infty}(0,T;\Ws) ,\\
        \mu \in L^{\infty}(0,T;\Vs) \cap  L^2(0,T;\Ws).
    \end{gather}
    In particular, $(\theta, \phi, \mu, \sigma)$ is a strong solution, i.e., it satisfies the system \eqref{eq:problem} a.e. in $Q$. Moreover, there exists a constant $C_2 > 0$  that depends on $M$ and on the other problem's data, such that the following estimate holds
    \begin{equation}
    \begin{split}
    &\norm{\theta}_{H^1(\Hs)\cap L^\infty(\Vs) \cap L^2(\Ws) \cap L^{\infty}(Q)} + \norm{\phi}_{W^{1,\infty}(\Vs^*)\cap H^1(\Vs) \cap  L^{\infty}(\Ws)}\\
    &\quad+ \norm{\mu}_{L^{\infty}(\Vs) \cap  L^2(\Ws)} \leq C_2.
    \end{split}
    \label{pier21}
    \end{equation}
\end{thm}

\begin{remark}
We point out that \eqref{hyp:reg_initial_condition} in particular yields \eqref{hyp:initial_condition} and $\hbeta(\phi_0) \in L^\infty (\Omega)$. Moreover, it is worth noting that, by standard Sobolev embedding results in dimensions $d =2,3$, the order parameter $\phi$ belongs to $L^{\infty}(Q)$.
\end{remark}

\noindent The $L^\infty$-bound for the component $\phi$ of the solution, combined with the local Lipschitz continuity of $\beta$, plays a crucial role in establishing a continuous dependence result, which, in particular, ensures the uniqueness property stated in the following theorem.

\begin{thm}[Continuous dependence]\label{thm:continuous_dependence}
    Suppose that the hypotheses \ref{hyp:constants}--\ref{hyp:beta_pi} are fulfilled, and that the initial data comply with \eqref{hyp:reg_initial_condition}. Then, for every assigned functions $\{u_i\}_{i=1,2}$ satisfying \eqref{hyp:control}, and any pair of initial data $\{(\theta_{0,i},\varphi_{0,i},\sigma_{0,i})\}_{i=1,2}$ satisfying \eqref{hyp:reg_initial_condition}, if we denote by $\{(\theta_i, \phi_i, \mu_i, \sigma_i)\}_{i=1,2}$ two corresponding strong solutions to \eqref{eq:problem}--\eqref{eq:initial_datum}, the following continuous dependence inequality holds
    \begin{equation}
        \begin{split}
            &\norm{\theta_1-\theta_2}_{L^{\infty}(\Hs) \cap L^2(\Vs)} + \norm{\phi_1-\phi_2}_{L^{\infty}(\Hs)\cap L^2(\Ws)}\\
            &\quad \quad  + \norm{\tau^{1/2}(\phi_1-\phi_2)}_{L^{\infty}(\Vs)}+ \norm{\sigma_1-\sigma_2}_{L^{\infty}(\Hs)\cap L^2(\Vs)}\\
            &\quad \leq C_3 \Big(\norm{\theta_{0,1} - \theta_{0,2}}_{\Hs} + \norm{\phi_{0,1} - \phi_{0,2}}_{\Hs} + \norm{\tau^{1/2}(\phi_{0,1} - \phi_{0,2})}_{\Vs}\\
            &\quad \quad + \norm{\sigma_{0,1}-\sigma_{0,2}}_{\Hs} + \norm{u_1-u_2}_{L^2(\Hs)}\Big)
        \end{split}
        \label{pier22}
    \end{equation}
    for a positive constant $C_3$ that depends on $M$ and on the other problem's data. Consequently, the strong solution found in~\Cref{thm:regularity} is unique.
\end{thm}

 \section{Existence of weak solutions}\label{sect:existence}
 The existence of weak solutions is proved through two levels of approximation. First, we introduce an approximated problem in which $\beta$ is replaced by a suitable Lipschitz continuous function $\betae$. Second, to show the existence of weak solutions of the approximate PDE system,  we employ a Faedo--Galerkin scheme. Finally, thanks to some a priori estimates, we separately pass to the limit, first in the Faedo--Galerkin scheme, and then as $\varepsilon \to 0$, finding a solution to the original problem. Notice that we need to introduce $\betae$ even if $\beta$ is smooth due to the wide class of possible growth that we are allowing: without it, we would not be able to pass to the limit in the corresponding term of the discretized system.

\subsection{The approximated problem}
 For $\varepsilon \in (0,1)$, we define the Moreau--Yosida approximation of $\hbeta$ as 
    \begin{equation}
        \hbetae(r) \coloneqq \min_{s \in \RR} \bigg\{\frac{1}{2\varepsilon}|s-r|^2+\widehat{\beta}(s )\bigg\} \qquad \forall r \in \RR,
        \label{pier-def-hbetae}
    \end{equation}
   and the Yosida regularization of $\beta$ as 
    \begin{equation}
        \betae \coloneqq \bigl(\hbetae\bigr)'.
        \label{pier-def-betae}    
        \end{equation}
They enjoy the following properties:
\begin{enumerate}[(i)]
    \item $\hbeta_{\varepsilon}$ is a $C^1(\RR)$ convex function with    \label{en:yosida_property_hat}    \begin{equation}\label{eq:yosida_property_hat}
        0 \leq \hbeta_{\varepsilon}(r) \leq \hbeta(r)
        \qquad \text{for all $r \in \RR$,}
    \end{equation}
    \item $\betae$ is monotone increasing and Lipschitz continuous with $\betae(0)=0$ and with Lipschitz constant bounded by $\varepsilon^{-1}$, \label{en:yosida_property}
    \item they satisfy the inequality \eqref{hyp:growth_beta} with the same constant $C_{\beta}$, i.e.,  \label{en:growth_betae}
    \begin{equation} \label{eq:growth_betae}
        |\betae(r)| \leq C_{\beta}(\hbetae(r) + 1)
        \quad \text{ for all $r \in \RR$.}
    \end{equation}
\end{enumerate}
In view of \ref{hyp:beta_pi},
the properties \ref{en:yosida_property_hat}--\ref{en:yosida_property} are well known and can be 
found in \cite[][Proposition~2.6, p.~28 and Proposition~2.11, p.~39]{brezis1973}. To prove 
\ref{en:growth_betae}, we need to introduce the resolvent operator associated with $\beta$, 
defined as
\begin{equation}
    J_{\varepsilon}(r)= ({\rm{Id}}+ \varepsilon \beta)^{-1}(r),
    \label{pier-risolvente}
\end{equation}
i.e., $J_{\varepsilon}(r)$ satisfy
\(
    J_{\varepsilon}(r) + \varepsilon \beta(J_{\varepsilon}(r)) = r
\)
for all $r \in \RR$. It can be proved (see \cite[][Proposition~2.11, p.~39]{brezis1973}) that $J_{\varepsilon} = \rm{Id} - \varepsilon \betae$.
It trivially follows that
\begin{equation}
    \beta(J_{\varepsilon}(r)) = \frac{r - J_{\varepsilon}(r)}{\varepsilon} = \betae(r)  \quad \text{ for all $r \in \RR$.}
    \label{pier-prop:betae}
    \end{equation}
Moreover, again by \cite[][Proposition~2.11, p.~39]{brezis1973}, we know that
\begin{equation*}
    \hbetae(r) = \frac{\varepsilon}{2} |\betae(r)|^2 + \hbeta(J_{\varepsilon}(r)),
\end{equation*}
thus $\hbetae(r) \geq \hbeta(J_{\varepsilon}(r))$. Putting these elements together and employing the inequality \eqref{hyp:growth_beta}, we obtain
\begin{equation*}
    |\beta_\varepsilon (r)| = |\beta (J_\varepsilon(r))| \leq C_{\beta} \left(1+ \hbeta (J_\varepsilon(r)) \right) \leq C_{\beta} \left(1+ \hbetae (r) \right),
\end{equation*}
so \ref{en:growth_betae} is proved.\\

\noindent The approximate problem is obtained from \eqref{eq:problem}--\eqref{eq:initial_datum} replacing $\beta$ with $\betae$ and its weak solutions are defined consequently as in \Cref{defn:weak_sol}.

\subsection{Existence of weak solutions for the approximated problem}
As already anticipated, we prove the existence of weak solutions of the approximated system through a Faedo--Galerkin space discretization.\\

\noindent \textbf{Faedo--Galerkin discretization.} We introduce the nondecreasing sequence of eigenvalues $\{\gamma_j\}_{j \in \NN_0}$ of the Laplace operator with homogeneous Neumann boundary conditions 
and the associated sequence of eigenvectors $\{e_j\}_{j \in \NN_0}$, which is a complete orthonormal system in $\Hs$, orthogonal in $\Vs$ and $\Ws$. Explicitly, for every $j \in \NN_0$,
\begin{equation*}
    \begin{cases}
        - \Delta e_j = \gamma_j e_j & \text{in } \Omega,\\
        \partial_{\nnu} e_j = 0  & \text{on } \partial \Omega,
    \end{cases}
\end{equation*}
with the additional properties 
\begin{equation*}
    \intprod{e_i}{e_j} = \delta_{ij}
    :=\begin{cases}1 & \hbox{if }\, i=j\\
       0 &\hbox{if }\, i\not=j 
       \end{cases}
    , \qquad  \intprod{\nabla e_i}{\nabla e_j} = \gamma_i \delta_{ij},
\end{equation*}
 for every $i, j \in \NN_0$. Recall that, by standard elliptic regularity results, $\{e_j\}_{j \in \NN_0}$ are smooth functions. We define
\begin{equation*}
    \Vs^n \coloneqq \text{span}\{e_0, \dots, e_n\} 
\end{equation*}
for every $n \in \NN_0$. Then,  $\{\Vs^n\}_{n \in \NN_0}$ is a nondecreasing sequence of subspaces, whose union is dense in $\Vs$ as well as in $\Hs$. With our notation, $\gamma_0 = 0$, $e_0 = |\Omega|^{-1/2}$, and $\Vs^0$ is the space of constant functions. We also introduce the projection of $\Hs$ onto $\Vs^n$ as 
\begin{equation*}
    P^n(v) \coloneqq \sum_{j = 0}^n \intprod{v}{e_j} e_j
\end{equation*}
for all $v \in \Hs$. Notice that there exists a constant $C>0$ such that
\begin{equation}
   \norm{P^n(v)}_X \leq C \norm{v}_X \quad \hbox{ for all $v \in X$, where $X = \Hs, \Vs, \Ws$.}
\label{pier-proiezione}
\end{equation}

\noindent For every $n \in \NN_0$, we aim at finding a quadruplet $(\thetan, \phin, \mun, \sigman)$ of the form
\begin{equation*}
    \begin{split}
        \thetan(x,t) = \sum_{j = 1}^n \thetanj(t)e_j(x), \quad \phin(x,t) = \sum_{j = 1}^n \phinj(t)e_j(x), \\
        \mun(x,t) = \sum_{j = 1}^n \munj(t)e_j(x), \quad \sigman(x,t) = \sum_{j = 1}^n \sigmanj(t)e_j(x)
    \end{split}
\end{equation*}
such that
\begin{subequations}\label{eq:gal_problem_eq_with_spaces}
    \begin{align}
        &\duality{\partial_t(\thetan + \ell\phin)}{v}_{\Vs} + \into  \nabla \thetan \cdot \nabla v \dx = \into u v \dx, \label{eq:gal_temperature_def}\\
        &\duality{\partial_t \phin}{v}_{\Vs} + \into  \nabla \mun \cdot \nabla v \dx = \into (\lambdap \sigman - \lambdaa - \lambdae \thetan) \h(\phin) v \dx, \label{eq:gal_ch1_def}\\
        &\begin{aligned}
             &\duality{\tau \partial_t \phin}{v}_{\Vs} + \into \nabla \phin\cdot \nabla v \dx + \into \left(\betae(\phin) + \pi(\phin) - \chi \sigman - \Lambda\thetan\right) v \dx\\
             &\quad = \into \mun v \dx,
        \end{aligned} \label{eq:gal_ch2_def}\\
        &\begin{aligned}
            &\duality{\partial_t\sigman}{v}_{\Vs} + \into \nabla \sigman \cdot \nabla v \dx - \chi \into \nabla \phin \cdot \nabla v \dx\\
            & \quad = \into \left(-\lambdac \sigman \h(\phin) + \lambdab (\sigma_B - \sigman) - \lambdad \sigman \k(\thetan)\right) v \dx,
        \end{aligned} \label{eq:gal_nutrient_def}
    \end{align}
\end{subequations}
 a.e. $t \in (0,T)$, for all $v \in \Vs^n$, and that fulfill the initial conditions
 \begin{equation}\label{eq:gal_initial_cond}
     \thetan(0) = \theta_{0,n} \coloneqq P^n(\theta_0), \quad \phin(0) = \phi_{0,n} \coloneqq P^n(\phi_0), \quad \sigman(0) = \sigma_{0,n} \coloneqq P^n(\sigma_0).
 \end{equation}
 Notice that, owing to the properties of the projection $P^n$, the following inequalities hold
 \begin{equation}
    \norm{\theta_{0,n}}_{\Hs} \leq C\norm{\theta_{0}}_{\Hs}, \quad \norm{\phi_{0,n}}_{\Vs} \leq C\norm{\phi_{0}}_{\Vs}, \quad \norm{\sigma_{0,n}}_{\Vs} \leq C\norm{\sigma_{0}}_{\Vs},
    \label{pier-ini-data}
 \end{equation}
 as well as the convergences
    \begin{alignat}{2}
    \theta_{0,n} \to \theta_0 &\quad \text{strongly} &&\quad \text{in } \Hs,\label{eq:gal_limit_initial_datum_theta}\\
    \phi_{0,n} \to \phi_0 &\quad \text{strongly} &&\quad\text{in } \Vs, \label{eq:gal_limit_initial_datum_phi}\\
    \sigma_{0,n} \to \sigma_0  &\quad\text{strongly} &&\quad\text{in } \Vs\label{eq:gal_limit_initial_datum_sigma}.
\end{alignat}
Here, we neglect the dependence on $\varepsilon$ in the notation $(\thetan, \phin, \mun, \sigman)$ on purpose, for brevity.\\

\noindent \textbf{Local-in-time existence.} We test each equation of the system \eqref{eq:gal_problem_eq_with_spaces}--\eqref{eq:gal_initial_cond} by $e_j$ for $j=0,\dots,n$, obtaining the following ODEs system
\begin{subequations}\label{eq:ode_system}
    \begin{align}
        &\difft (\thetanj +\ell\phinj) + \gamma_j \thetanj = \intprod{u}{e_j},\label{eq:ode_system_temp}\\
        &\difft \phinj + \gamma_j \munj = \intprod{(\lambdap \sigman - \lambdaa - \lambdae \thetan)\h(\phin)}{e_j},\label{eq:ode_system_ch1}\\
        &\tau \difft \phinj + \gamma_j \phinj + \intprod{\betae(\phin) + \pi(\phin)}{e_j} -\chi \sigmanj -\Lambda \thetanj
        = \munj,\label{eq:ode_system_ch2}\\
        &\difft \sigmanj + (\gamma_j + \lambdab) \sigmanj
        - \chi \gamma_j \phinj
        = \intprod{\left(-\lambdac \sigman \h(\phin) + \lambdab \sigma_B - \lambdad \sigman \k(\thetan)\right)}{e_j} \label{eq:ode_system_nutrient},\\
        &\thetanj(0) = \intprod{\theta_{0,n}}{e_j}, \quad \phin(0) = \intprod{\phi_{0,n}}{e_j},\quad \sigmanj(0) = \intprod{\sigma_{0,n}}{e_j}.
    \end{align}
\end{subequations}
Notice that $\munj$ is an auxiliary variable and can be removed from the system by substituting $\munj$, whose expression is given by equation \eqref{eq:ode_system_ch2}, into equation \eqref{eq:ode_system_ch1}. Then we can recover the expression for $\difft \phinj$ and replace it in equation \eqref{eq:ode_system_temp}. Hence, we obtain a $3n$-equations first-order ODEs system in the variables $\thetanj$, $\phinj$, and $\sigmanj$ for $j=1,\dots,n$ with locally Lipschitz continuous nonlinearities. Thus, the existence of $H^1(0,T^n)$ solutions follows from Carath\'eodory existence theorem, for a certain $T^n \leq T$. The solution is not global because the nonlinearities are not globally Lipschitz continuous.  Then, we retrieve $\munj \in L^2(0,T^n)$ by equation \eqref{eq:ode_system_ch2}. 
In the following, we will derive some a priori estimates independent of $n$ that will allow us to extend the solution to all $[0,T]$, and, at the same time, to recover enough compactness to pass to the limit as $n$ goes to infinity. Keeping in mind that our final goal is passing to the limit as $\varepsilon \to 0$, we will also be careful in tracking the dependence on $\varepsilon$.\\

\noindent \textbf{First a priori estimate.} We integrate in time 
equation~\eqref{eq:gal_temperature_def} for a fixed $v \in \Vs^n$, obtaining the equation
\begin{equation}\label{eq:eq:gal_temperature_weaker_def}
    \into (\thetan +\ell \phin) v \dx + \into \nabla(1*_t\thetan)\cdot \nabla v \dx = \into (1*_t u)v + \into (\theta_{0,n} +\ell \phi_{0,n})v \dx.
\end{equation}
Here, the notation $*_t$ denotes the convolution with respect to time, i.e., 
\begin{equation*}
    1*_t\thetan \coloneqq \intt \thetan(\cdot,s) \ds,\quad 1*_t u \coloneqq \intt u(\cdot,s) \ds.
\end{equation*}
We take $v = \thetan$ in equation \eqref{eq:eq:gal_temperature_weaker_def} and multiply each side of the equality by a constant $R>0$, fixed but yet to be determined. We choose $v = \phin$ in \eqref{eq:gal_ch1_def}, $v = - \Delta \phin$ in \eqref{eq:gal_ch2_def}, and $v = \sigman$ in \eqref{eq:gal_nutrient_def}. We add the resulting equalities and, after some simplification, we obtain:
\begin{align*}
    &R \into |\thetan|^2 \dx + \frac{R}{2} \difft \into |\nabla (1 *_t \thetan)|^2 \dx + \frac{1}{2}\difft\into |\phin|^2 \dx + \frac{\tau}{2} \difft \into |\nabla \phin|^2 \dx\\
    &\quad\quad+ \into |-\Delta \phin|^2 \dx + \into \betae'(\phin)|\nabla \phin|^2 \dx + \frac{1}{2} \difft \into |\sigman|^2 \dx + \into |\nabla \sigman|^2 \dx\\
    &\quad = R \into (1*_t u)\thetan \dx + R \into (\theta_{0,n}+ \ell\phi_{0,n})\thetan \dx - R \ell \into \phin \thetan \dx\\
    &\quad \quad + \into (\lambdap\sigman- \lambdaa - \lambdae\thetan)\h(\phin)\phin \dx
   -\into \pi(\phin)(-\Delta \phin)\dx
    + 2 \into \chi \sigman(-\Delta \phin) \dx\\
    &\quad \quad+ \Lambda\into \thetan(-\Delta \phin) \dx - \into \left(\lambdac \h(\phin) + \lambdad \k(\thetan)\right)|\sigman|^2\dx + \into  \lambdab (\sigma_B -\sigman) \sigman \dx.
\end{align*}
\noindent After applying the H\"older and the Young inequalities to the right-hand side, and rearranging some terms, we end up with
\begin{align*}
        &R \into |\thetan|^2 \dx  + \into |-\Delta \phin|^2 \dx + \into \betae'(\phin)|\nabla \phin|^2 \dx + \into |\nabla \sigman|^2 \dx
        + \into  \lambdab |\sigman|^2 \dx\\
        &\quad \quad+ \frac{1}{2} \difft \left( R\into |\nabla (1 *_t \thetan)|^2 \dx + \into |\phin|^2 \dx + \tau \into |\nabla \phin|^2 \dx  + \into |\sigman|^2 \dx \right)\\
        &\quad = R \into \left(1*_t u + \theta_{0,n} + \ell \phi_{0,n} - \ell \phin \right)\thetan \dx - \into \lambdae\h(\phin) \phin \thetan \dx\\
        &\quad \quad + \into (\lambdap \sigman - \lambdaa)\h(\phin)\phin \dx + \into \left(- \pi(\phin) + 2 \chi \sigman + \Lambda \thetan\right)(-\Delta \phin) \dx\\
        &\quad \quad  - \into \left(\lambdac \h(\phin) + \lambdad \k(\thetan)\right)|\sigman|^2\dx + \into \lambdab \,\sigma_B \,\sigman\dx\\
        & \quad \leq \bigg(\frac{R}{2} + \frac{1}{2} + \frac{\Lambda^2}{2} \bigg) \into |\thetan|^2 \dx + \frac{1}{2} \into |-\Delta \phin|^2 \dx \\
        & \quad \quad + C_R \bigg( \into |\theta_{0,n}|^2 \dx + \into |\phi_{0,n}|^2 \dx + \into |1*_t u|^2 \dx + \into |\phin|^2 \dx \bigg)\\
        & \quad \quad  + C\bigg( \into |\sigma_B|^2 \dx + \into |\sigman|^2 \dx +  1  \bigg), 
\end{align*}
where we have exploited the boundedness of $\h$ and $\k$ given by hypothesis \ref{hyp:given_function} and the Lipschitz continuity of $\pi$ (cf.~\eqref{hyp:pi_lip}). Now we choose $R = R(\Lambda)$ such that $R > 1 + \Lambda^2$, move the first two addends of the right-hand side to the left-hand side, and integrate in time over the interval $(0,t)$ for $t \leq T^n \leq T$. Moreover, we recall that $\betae'$ is nonnegative because $\betae$ is monotone. Thus, the corresponding integral is also nonnegative, and we can get rid of it. After renaming the constants, we have:
\begin{equation*}
   	\begin{split}
            &\intto |\thetan|^2 \dx \ds + \intto |-\Delta \phin|^2 \dx \ds + \intto |\nabla \sigman|^2 \dx\ds + \intto |\sigman|^2 \dx\ds\\
            &\quad \quad + \into |\nabla(1*_t\thetan)|^2 \dx + \into |\phin|^2 \dx + \tau \into |\nabla \phin|^2 \dx + \into |\sigman|^2 \dx\\
            &\quad \leq C_T \bigg( 1 + \into \left(|\theta_{0,n}|^2 + |\phi_{0,n}|^2 + \tau|\nabla \phi_{0,n}|^2 + |\sigma_{0,n}|^2 \right) \dx\\
            &\quad \quad + \intTo \left(|1*_t u|^2 + |\sigma_B|^2 \right)\dx \dt + \intto \left(|\phin|^2 + |\sigman|^2\right) \dx \ds \bigg).
   	\end{split}
\end{equation*}
Applying the Gronwall inequality leads to
\begin{equation}\label{eq:gal_1st_est}
   	\begin{split}
   		&\norm{\thetan}_{L^2(0,T^n;\Hs)} + \norm{1*_t\thetan}_{L^{\infty}(0,T^n;\Vs)} + \norm{\phin}_{L^{\infty}(0,T^n;\Hs) \cap L^2(0,T^n;\Ws)}\\
   		&\quad + \norm{\tau^{1/2} \phin}_{L^{\infty}(0,T^n;\Vs)} + \norm{\sigman}_{L^{\infty}(0,T^n;\Hs) \cap L^2(0,T^n;\Vs)} \leq C,
   	\end{split}
\end{equation}
where $C$ depends on the assigned data of the problem (in particular, on $\norm{1 *_t u}_{L^2(\Hs)}$), but not on $T^n$ and $\varepsilon$.\\ 

\noindent \textbf{First a priori estimate consequences.} From estimate \eqref{eq:gal_1st_est}, it trivially follows that
\begin{equation}\label{eq:gal_phi_mean_value_est}
   	\norm{\ophin}_{L^{\infty}(0,T^n)} \leq C.
\end{equation}
Moreover, testing the equation \eqref{eq:gal_ch1_def} with $v = |\Omega|^{-1}$ and estimating the right-hand side with its $\Hs$-norm, we obtain
\begin{equation}\label{eq:gal_dt_phi_mean_value_est}
   	\norm{\overline{\partial_t\phin}}_{L^2(0,T^n)} \leq C \norm{(\lambdap \sigman - \lambdaa - \lambdae \thetan)\h(\phin)}_{L^2(0,T^n;\Hs)} \leq C. 
\end{equation} 

 \noindent \textbf{Second a priori estimate.}
 We preliminarily notice that if $v \in \Vs^n$ and has zero mean value, then $\N v \in \Vs^n$. Indeed, since $v$ belongs to $\Hs$, then $\N v$ belongs to $\Ws$ and the following equality obviously holds
 \begin{equation*}
     -\Delta(\N v) = \sum_{j = 0}^{+\infty} \intprod{-\Delta(\N v)}{e_j} e_j = \sum_{j = 0}^{+\infty} \gamma_j \intprod{ \N v}{e_j} e_j.
 \end{equation*}
 On the other hand, we have that
 \begin{equation*}
     - \Delta (\N v) = v = \sum_{j=0}^n \intprod{v}{e_j} e_j.
 \end{equation*}
 Consequently, since these two expressions must coincide and $\gamma_j > 0$ for every $j \in \NN$, we deduce that $\intprod{\N v}{e_j}$ is equal to zero for every $j \geq n + 1$. Thus, $\N v \in \Vs^n$.\\
Now, we take the difference of equation~\eqref{eq:gal_ch1_def} with its mean value and test it with $\N (\phin - \ophin)$, then we add the resultant to equation \eqref{eq:gal_ch2_def} tested with $(\phin - \ophin)$. In view of the property~\eqref{enum:eqN}, we note the cancellation of the terms involving $\mun$ as well as of the scalar products of a mean value and of $\N (\phin - \ophin)$ or 
$(\phin - \ophin)$. Thus, we obtain
\begin{equation*}
    \begin{split}
        &\frac{1}{2}  \difft \left(\into |\nabla \N(\phin - \ophin)|^2 \dx  + \tau \into |\phin - \ophin|^2 \dx \right)\\
        &\quad \quad  + \into |\nabla \phin|^2 \dx + \into \betae(\phin)(\phin - \ophin) \dx \\
        &\quad =  \into (\lambdap \sigman - \lambdaa - \lambdae \thetan)\h(\phin)\N(\phin - \ophin) \dx\\
        &\quad \quad- \into \pi(\phin) (\phin - \ophin) \dx + \chi \into \sigman (\phin - \ophin)\dx + \Lambda \into \thetan (\phin - \ophin).
    \end{split}
\end{equation*}     
Regarding the left-hand side, we observe that
\begin{equation*}
   	\into \betae(\phin)(\phin-\ophin) \dx \geq \into \hbetae(\phin) \dx - \into \hbetae(\ophin) \dx,	
\end{equation*}
due to the fact that $\betae$ is the subdifferential of $\hbetae$. Then, we integrate in time over the interval $(0,t)$ for $t \leq T^n \leq T$, and estimate the terms on the right-hand side using the H\"older and the Young inequalities, obtaining
\begin{equation*}
   	\begin{split}
 			&\frac{1}{2} \left(\into |\nabla \N(\phin - \ophin)|^2 \dx  + \tau \into |\phin - \ophin|^2 \dx \right)\\
 			&\quad \quad  + \intto |\nabla \phin|^2 \dx \ds + \intto \hbetae(\phin) \dx \ds \\
 			&\quad \leq  \intto \hbetae(\ophin) \dx \ds+ C \intto |\N(\phin-\ophin)|^2 \dx \ds\\
 			&\quad \quad + C \bigg(\intto |\sigman|^2 + |\thetan|^2 + |\phin|^2 + |\phin -\ophin|^2 \dx \ds +1 \bigg) \eqqcolon I_1 + I_2 + I_3, 
   	\end{split}
\end{equation*}
where we used the boundedness of $\h$ by hypothesis \ref{hyp:given_function}, and the fact that $\pi$ is Lipschitz continuous by hypothesis \ref{hyp:beta_pi}. To handle $I_1$, we recall that $\hbetae \leq \hbeta$, where $\hbeta$ is continuous by hypothesis \ref{hyp:beta_pi}, and that, by \eqref{eq:gal_phi_mean_value_est}, $\ophin$ is uniformly bounded. Therefore, it follows that~$I_1 \leq C$. We use the Poincaré inequality to treat $I_2$, taking into account that $\N(\phin - \ophin)$ has zero mean value. Finally, $I_3$ can be uniformly bounded thanks to \eqref{eq:gal_1st_est} and \eqref{eq:gal_phi_mean_value_est}.  Thus, we have
\begin{equation*}
    \begin{split}
        &\into |\nabla \N(\phin - \ophin)|^2 \dx  + \tau \into |\phin - \ophin|^2 \dx + \intto |\nabla \phin|^2 \dx \ds + \intto \hbetae(\phin) \dx \ds \\
        &\quad \leq  C \left(\intto |\nabla \N(\phin-\ophin)|^2 \dx \dt + 1 \right)
    \end{split}
\end{equation*}
from which, through the Gronwall Lemma, we obtain
\begin{equation}\label{eq:gal_2nd_est}
  	 \norm{\hbetae(\phin)}_{L^1(\Omega \times (0,T^n))} \leq C.
\end{equation} 

\noindent \textbf{Second a priori estimate consequences.} From \eqref{eq:gal_2nd_est}, recalling the $\hbetae$ growth 
property~\eqref{eq:growth_betae}, it trivially follows that
\begin{equation}\label{eq:gal_betae_est}
   	\norm{\betae(\phin)}_{L^1(\Omega \times (0,T^n))} \leq C_{\beta}(1 + \norm{\hbetae(\phin)}_{L^1(\Omega \times (0,T^n))}) \leq C.
\end{equation}
Consequentially, taking $v = |\Omega|^{-1}$ in equation \eqref{eq:gal_ch2_def}, recalling that $\overline{\partial_t \phin}$ is estimated from \eqref{eq:gal_dt_phi_mean_value_est} and thanks to the other estimates from \eqref{eq:gal_1st_est}, we deduce that
\begin{equation}\label{eq:gal_mu_mean_value_est_pt1}
   	\norm{\omun}_{L^1(0,T^n)} \leq C.
\end{equation}

\noindent \textbf{Third a priori estimate.} We take $v =(\Lambda/\ell) \thetan$ in equation \eqref{eq:gal_temperature_def}, $v = \mun$ in \eqref{eq:gal_ch1_def}, $v = \partial_t \phin$ in \eqref{eq:gal_ch2_def}, and $v = \partial_t \sigman$ in \eqref{eq:gal_nutrient_def}. We sum the resulting equalities, noting a cancellation, and integrate in time over $(0,t)$ for a $t \leq T^n \leq T$, obtaining:
\begin{align*}\label{eq:gal_est1_testing}
    &\frac{\Lambda}{2\ell} \into |\thetan|^2 \dx + \frac{1}{2}\into |\nabla \phin|^2 \dx + \into \hbetae(\phin) \dx \\
    &\quad \quad + \frac{1}{2}\into |\nabla \sigman|^2 \dx + \frac{\Lambda}{\ell}\intto |\nabla \thetan|^2 \dx \ds + \intto |\nabla \mun|^2 \dx \ds \\
    &\quad \quad+ \tau \intto |\partial_t \phin|^2 \dx \ds + \intto |\partial_t \sigman|^2 \dx \ds\\
    &\quad = \frac{\Lambda}{2\ell}
    \into |\theta_{n,0}|^2 \dx + \frac{1}{2}\into |\nabla \phi_{n,0}|^2 \dx 
    + \into \hbetae(\phi_{n,0}) \dx\\
    &\quad \quad + \frac{1}{2}\into |\nabla \sigma_{n,0}|^2 \dx + \into \hpi(\phi_{n,0}) - \hpi(\phin)\dx
    +\frac{\Lambda}{\ell}\intto u \thetan \dx \ds  \\
    &\quad \quad + \intto (\lambdap \sigman - \lambdaa - \lambdae \thetan) \h(\phin)\mun \dx \ds + \intto \chi \sigman \partial_t \phin \dx \ds\\
    &\quad \quad- \intto \lambdac \sigman \h(\phin)\partial_t \sigman \dx \ds + \intto \lambdab(\sigma_B - \sigman) \partial_t \sigman \dx \ds\\
    &\quad \quad- \intto \lambdad \sigman \k(\thetan) \partial_t \sigman \dx \ds + \intto \chi (-\Delta \phin) \partial_t \sigman \dx \ds.
\end{align*}
We aim to bound from above the right-hand side of this equality. First, we deal with the terms related to the initial data of the Galerkin discretized system. To do so, we recall that since they are projections of the original initial data, their $\Hs$ and $\Vs$ norms can be estimated with a constant independent of $n$. We notice that $\hpi$ has at most quadratic growth thanks to hypothesis \ref{hyp:beta_pi}. Moreover, the Moreau--Yosida approximation $\hbetae$ satisfies 
\begin{equation}
    \hbetae(r) \leq \hbetae(0) + \betae(0)r +
    \frac1{2\varepsilon}r^2 \leq \frac1{2\varepsilon}r^2
    \quad\text{for every }r\in \mathbb{R},
\end{equation}
because of the definitions~\eqref{pier-def-hbetae}, \eqref{pier-def-betae} and the related property \ref{en:yosida_property} (cf. also~\ref{hyp:beta_pi}), so here---just in order to bound the term $$\into \hbetae(\phi_{n,0}) \dx$$---we obtain an estimate which is not independent of~$\varepsilon$. 
Most of the other terms on the right-hand side can be estimated simply through the H\"older and the Young inequalities, and the first a priori 
estimate~\eqref{eq:gal_1st_est}, leading to
\begin{align}
            &\frac{\Lambda}{2\ell} \into |\thetan|^2 \dx + \frac{1}{2}\into |\nabla \phin|^2 \dx + \into \hbetae(\phin) \dx \notag\\
            &\quad \quad+ \frac{1}{2}\into |\nabla \sigman|^2 \dx + \frac{\Lambda}
            {\ell}\intto |\nabla \thetan|^2 \dx \ds + \intto |\nabla \mun|^2 \dx \ds \notag\\
            &\quad \quad + \tau \intto |\partial_t \phin|^2 \dx \ds + \intto |\partial_t \sigman|^2 \dx \ds \notag\\
            &\quad \leq C_{0,\varepsilon} +C \into \bigl(|\phin|^2 +1 \bigr)\dx + \frac{\Lambda}{2\ell}\intto |u|^2 \dx \ds + \frac{\Lambda}{2\ell}\intto |\thetan|^2  \dx\ds \notag\\ 
            &\quad \quad  + \intto (\lambdap \sigman - \lambdaa - \lambdae \thetan) \h(\phin)\mun \dx \ds + \intto \chi \sigman \partial_t \phin \dx \ds \label{eq:gal_2nd_est_pt1}\\
            &\quad \quad + \frac{1}{4} \intto |\partial_t \sigman|^2 \dx \ds + C\intto |\sigman|^2 + |\sigma_B|^2 + |-\Delta \phin|^2\dx \ds \notag\\
            &\quad \leq C_{0,\varepsilon} + C + \frac{1}{4} \intto |\partial_t \sigman|^2 \dx \ds \notag\\
            &\quad \quad + \intto (\lambdap \sigman - \lambdaa - \lambdae \thetan) \h(\phin)\mun \dx \ds + \intto \chi \sigman \partial_t \phin \dx \ds \notag\\
            &\quad\eqcolon C_{0,\varepsilon} + C + \frac{1}{4} \intto |\partial_t \sigman|^2 \dx \ds + I_4 + I_5.\notag
\end{align}
At this point, we focus on the terms $I_4$ and $I_5$. To treat $I_4$, we employ the H\"older inequality and the fact that $\norm{\sigman}_{L^{\infty}(0,T^n;\Hs)}$ is uniformly bounded by \eqref{eq:gal_1st_est}. Then, we use the Poincaré, the H\"older, the Young inequalities and the previous estimates \eqref{eq:gal_1st_est}, \eqref{eq:gal_mu_mean_value_est_pt1}. We have:
\begin{equation}\label{eq:gal_2nd_est_pt2}
   	\begin{split}
   			I_4 & \leq C\intt ( \norm{\sigman}_{\Hs} + \norm{\thetan}_{\Hs} + 1) \norm{\mun}_{\Hs} \ds \leq C\intt ( \norm{\thetan}_{\Hs} + 1) \norm{\mun}_{\Hs} \ds\\
   			& \leq C\intt (\norm{\thetan}_{\Hs} + 1) \norm{\nabla \mun}_{\Hs} \ds + C\intt (\norm{\thetan}_{\Hs} + 1)|\omun| \ds \\
   			&\leq \frac{1}{2} \intt \norm{\nabla \mun}_{\Hs}^2 \ds + C + C \intt |\omun|\, \norm{\thetan}_{\Hs}^2 \ds.
   	\end{split}
\end{equation}
We handle the term $I_5$ integrating by parts with respect to time and then through the H\"older and Young inequalities. Explicitly, we have
\begin{equation}\label{eq:gal_2nd_est_pt3}
   	\begin{split}
   		I_5 &= - \intto \chi \partial_t \sigman \phin \dx \ds + \into \chi \sigman \phin \dx - \into \chi \sigma_{n,0}\phi_{n,0} \dx\\\
   		&\leq \frac{1}{4} \intto |\partial_t \sigman|^2 \dx \ds + C \left( \intto |\phin|^2 \dx \ds + \into |\sigman|^2  + |\phin|^2 \dx +1\right)\\
   		&\leq \frac{1}{4} \intto |\partial_t \sigman|^2 \dx \ds + C,
   	\end{split}
\end{equation}  
where in the last inequality we used the first a priori estimate \eqref{eq:gal_1st_est}. We return to \eqref{eq:gal_2nd_est_pt1}, and make use of \eqref{eq:gal_2nd_est_pt2} and \eqref{eq:gal_2nd_est_pt3}. Rearranging some terms and renaming the constants, we end up with
\begin{equation*}
   	\begin{split}
   			&\into |\thetan|^2 \dx + \into |\nabla \phin|^2 \dx + \into \hbetae(\phin) \dx + \into |\nabla \sigman|^2 \dx\\
   			&\quad \quad + \intto |\nabla \thetan|^2 \dx \ds + \intto |\nabla \mun|^2 \dx \ds \dx\\
   			&\quad \quad + \tau \intto |\partial_t \phin|^2 \dx \ds + \intto |\partial_t \sigman|^2 \dx \ds\\
   			&\quad \leq  C_{0,\varepsilon} + C  +  C \intt 
            |\omun|\, \norm{\thetan}_{\Hs}^2\ds.
   	\end{split}
\end{equation*} 
Thanks to the Gronwall inequality and the fact that $\norm{\omun}_{L^1(0,T^n)}$ is uniformly bounded from \eqref{eq:gal_mu_mean_value_est_pt1}, we deduce
\begin{equation} \label{eq:gal_3rd_est}
	\begin{split}
		&\norm{\thetan}_{L^{\infty}(0,T^n;\Hs) \cap L^2(0,T^n;\Vs)} + \norm{\phin}_{L^{\infty}(0,T^n;\Vs)}\\
		&\quad   + \norm{\tau^{1/2}\phin}_{H^1(0,T^n;\Hs)}+ \norm{\hbetae(\phin)}_{L^{\infty}(0,T^n;L^1(\Omega))}\\
		&\quad   + \norm{\nabla \mun}_{L^2(0,T^n;\Hs)} +  \norm{\sigman}_{H^1(0,T^n;\Hs) \cap L^{\infty}(0,T^n;\Vs)} \leq C_{\varepsilon}.
	\end{split}
\end{equation}

\noindent \textbf{Third a priori estimate consequences.}      
 Owing to the estimate \eqref{eq:gal_3rd_est} and
 exploiting the growth property~\eqref{eq:growth_betae}, it is straightforward to infer that 
\begin{equation}\label{eq:gal_beta_est_infty}
    \norm{\betae(\phin)}_{L^{\infty}(0,T^n;L^1(\Omega))} \leq C_{\varepsilon}.
\end{equation}
Thus, proceeding as before and choosing $v = |\Omega|^{-1}$ in the equation \eqref{eq:gal_ch2_def}, by comparing the terms we have that
\begin{equation}\label{eq:gal_mu_mean_value_est}
   	\norm{\omun}_{L^2(0,T^n)} \leq C_{\varepsilon},
\end{equation}
thanks to~\eqref{eq:gal_3rd_est}. Then, by the Poincar\'e inequality, it follows that 
\begin{equation}\label{eq:gal_mu_est}
   	\norm{\mun}_{L^2(0,T^n;\Vs)} \leq C_{\varepsilon}.
\end{equation}
Next, we want to show that
\begin{equation}\label{eq:gal_dt_phin_dualV}
   	\norm{\partial_t \phin}_{L^2(0,T^n;\Vs^*)} \leq C_{\varepsilon}.
\end{equation}
To do so, we consider a $v \in \Vs$ and proceed as follows
\begin{equation}
   	\begin{split}
   		&\left|\duality{\partial_t \phin}{v}_{\Vs} \right|= 	\left|\duality{\partial_t \phin}{P^n(v)}_{\Vs} \right|\\
   		&\quad= \left|-\into \nabla \mun \cdot \nabla \left[P^n(v)\right] \dx + \into (\lambdap \sigman - \lambdaa - \lambdae \thetan)\h(\phin)P^n(v) \dx \right|\\
   		&\quad\leq C \left(\norm{\nabla \mun}_{\Hs} + \norm{\sigman}_{\Hs} + \norm{\thetan}_{\Hs} + 1\right)\norm{P^n(v)}_{\Vs}\\
   		&\quad\leq C \left(\norm{\nabla \mun}_{\Hs} + \norm{\sigman}_{\Hs} + \norm{\thetan}_{\Hs} + 1\right)\norm{v}_{\Vs},
   	\end{split}
    \label{pier14}
\end{equation}
where, besides the usual calculations, we exploit the fact that $\partial_t \phin$ satisfies equation \eqref{eq:gal_ch1_def} and that, even if $v$ does not belong to the test functions space $V^n$, its projection does. From this inequality, we find that
\begin{equation*}
    \begin{split}
        \norm{\partial_t \phin}_{L^2(0,T^n;\Vs^*)}^2 &= \intTn \norm{\partial_t \phin}_{\Vs^*}^2 \dt \\
        &\leq C \intTn \left(\norm{\nabla \mun}_{\Hs}^2 + \norm{\sigman}_{\Hs}^2 + \norm{\thetan}_{\Hs}^2 + 1\right) \dt \leq C_{\varepsilon},
    \end{split}
\end{equation*}
so \eqref{eq:gal_dt_phin_dualV} follows.
Proceeding similarly, it is easy to prove that
\begin{equation}\label{eq:gal_dt_thetan_dualV}
  			\norm{\partial_t \thetan}_{L^2(0,T^n;\Vs^*)} \leq C_{\varepsilon}.
\end{equation}
Finally, we take $v = - \Delta \sigman$ in equation \eqref{eq:gal_nutrient_def} and integrate in time, obtaining
\begin{equation*}
  		\begin{split}
  			&\norm{-\Delta \sigman}_{L^2(0,T^n;\Hs)}^2 \dt \\
            &\quad = - \intTno \partial_t \sigman (-\Delta \sigman) \dx \dt + \chi \intTno (-\Delta \phin) (-\Delta \sigman) \dx \dt\\
  			&\quad \quad  + \intTno \left(-\lambdac \sigman\h(\phin) + \lambdab(\sigma_B - \sigman) - \lambdad \sigman \k(\thetan)\right)(-\Delta \sigman)\dx \dt\\
  			&\quad \leq C \left(\norm{\sigman}_{H^1(0,T^n;\Hs)} + \norm{-\Delta\phin}_{L^2(0,T^n;\Hs)} + \norm{\sigma_B}_{L^2(0,T^n;\Hs)}   \right)\norm{-\Delta \sigman}_{L^2(0,T^n;\Hs)}\\
  			&\quad \leq C_{\varepsilon} \norm{-\Delta \sigman}_{L^2(0,T^n;\Hs)},
  		\end{split}
\end{equation*}
where we used the H\"older inequality and the previously proved estimates \eqref{eq:gal_1st_est} and \eqref{eq:gal_3rd_est}. Thus, $\norm{-\Delta \sigman}_{L^2(0,T^n;\Hs)}$ is uniformly bounded and, by elliptic regularity, we have that
\begin{equation} \label{eq:gal_sigman_H2}
  		\norm{\sigman}_{L^2(0,T^n;\Ws)} \leq C_{\varepsilon}.
\end{equation} 

\noindent \textbf{Global-in-time existence.}
Let $\overline{T}^{n} \leq T$ denote the maximal existence time of the local-in-time solution $(\thetan, \phin, \mun, \sigman)$. Suppose, by contradiction, that $\overline{T}^{n} < T$. 
Thanks to the uniform (in $n$) a priori estimates established earlier and standard embedding results, by continuity the solution is defined at time $\overline{T}^{n}$ and satisfies the requirement on the initial values (cf.~\eqref{pier-ini-data} and~\eqref{hyp:initial_condition})
\begin{equation*}
 \norm{\theta^{n}\bigl(\overline{T}^{n}\bigr) }_{\Hs} + \norm{\phi^{n}\bigl(\overline{T}^{n}\bigr)}_{\Vs} + \norm{\sigma^{n}\bigl(\overline{T}^{n}\bigr)}_{\Vs} \leq C. 
\end{equation*}
Therefore, we can take the solution at $t = \overline{T}^{n}$ as a new initial datum for the system of differential equations and extend the solution beyond $\overline{T}^{n}$ by continuity.
This leads to a contradiction, since it implies that $\overline{T^n}$ is not maximal. We conclude that $\overline{T}^n = T$.\\

\noindent \textbf{Passage to the limit as $n \to \infty$.} In the previous steps, we proved that the solution of the Faedo--Galerkin discretized system satisfies the following bound
\begin{align*}
  	&\norm{\thetan}_{H^1(\Vs^*) \cap L^{\infty}(\Hs) \cap L^2(\Vs)} +
  	\norm{\phin}_{H^1(\Vs^*) \cap L^{\infty}(\Vs) \cap L^2(\Ws)}\\
  	&\quad+ \norm{\mun}_{L^2(\Vs)} +
  	\norm{\sigman}_{H^1(\Hs) \cap L^{\infty}(\Vs) \cap L^2(\Ws)} \leq C_{\varepsilon}.
\end{align*}
Thus, from standard compactness results (i.e., Banach--Alaoglu and Aubin--Lions theorems), we deduce that there exists a quadruplet $(\thetae,\phie, \mue, \sigmae)$ satisfying the regularity of \Cref{defn:weak_sol} such that, for $n \to +\infty$, along a nonrelabelled subsequence, the following convergences hold:
  	\begin{align}
  		&\thetan \to \thetae &&\text{weakly-}\ast &&\text{in } H^1(0,T;\Vs^*) \cap L^{\infty}(0,T;\Hs) \cap L^2(0,T;\Vs), \label{eq:limit_thetan_w}\\
  		& &&\text{strongly} &&\text{in } L^2(0,T;\Hs), \label{eq:limit_thetan_s}\\
  		& &&\text{a.e.} &&\text{in } Q, \label{eq:limit_thetan_ae}\\
  		&\phin \to \phie &&\text{weakly-}\ast &&\text{in } H^1(0,T;\Vs^*) \cap L^{\infty}(0,T;\Vs) \cap L^2(0,T;\Ws), \label{eq:limit_phin_w}\\
  		& &&\text{strongly} &&\text{in } C^0([0,T]; \Hs) \cap L^2(0,T;\Vs), \label{eq:limit_phin_s}\\
  		& &&\text{a.e.} &&\text{in }  Q, \label{eq:limit_phin_ae}\\
  		&\mun \to \mue &&\text{weakly} &&\text{in } L^2(0,T;\Vs), \label{eq:limit_mun_w}\\
  		&\sigman \to \sigmae &&\text{weakly-}\ast &&\text{in } H^1(0,T;\Hs) \cap L^{\infty}(0,T;\Vs) \cap L^2(0,T;\Ws), \label{eq:limit_sigman_w}\\
  		& &&\text{strongly} &&\text{in } C^0([0,T];\Hs) \cap L^2(0,T;\Vs), \label{eq:limit_sigman_s}\\
  		& &&\text{a.e.} &&\text{in } Q. \label{eq:limit_sigman_ae}
  	\end{align}
We integrate in time the equations \eqref{eq:gal_temperature_def}--\eqref{eq:gal_nutrient_def} over the interval $(0,T)$, after choosing as a test function $P^n(v) \in L^2(0,T;\Vs^n)$  for a generic fixed $v \in L^2(0,T;\Vs)$. We recall that 
\begin{equation}\label{eq:strong_conv_test_ft}
    P^n(v) \to v \quad \text{strongly} \quad \text{in } L^2(0,T;\Vs). 
\end{equation} 
All the linear terms pass to the limit thanks to the weak convergences \eqref{eq:limit_thetan_w}, \eqref{eq:limit_phin_w}, \eqref{eq:limit_mun_w}, \eqref{eq:limit_sigman_w}, and to the strong convergence \eqref{eq:strong_conv_test_ft}. Let's discuss only the nonlinear ones. In the mass source term in the equation \eqref{eq:gal_ch1_def}, we have that
\begin{equation*}
    (\lambdap \sigman - \lambdaa - \lambdae \thetan)\h(\phin) \quad \to \quad (\lambdap \sigmae - \lambdaa - \lambdae \thetae)\h(\phie)
\end{equation*} 
weakly in $L^2(0,T;\Hs)$. In fact, the term between parentheses strongly converges in $L^2(0,T;\Hs)$ and a.e. in $Q$; besides, $\h$ is continuous and bounded by hypothesis \ref{hyp:nonlinearities} and $\phin \to \phi$ a.e. in $Q$.
  Thus, we have weak convergence in $L^2(0,T;\Hs)$ by uniform boundedness, and a.e. convergence, with identification of the limit (see, e.g., \cite[Lemme~1.3, p.~12]{Lions_69}).
 Regarding the equation \eqref{eq:gal_ch2_def}, we notice that $\betae + \pi$ is Lipschitz continuous. Thus, the strong convergence \eqref{eq:limit_phin_s} is enough to pass to the limit. Finally, the term
  	\begin{equation*}
  		-\lambdac \sigman \h(\phin) - \lambdad \sigman \k(\thetan) \quad  \to \quad -\lambdac \sigmae  \h(\phie) - \lambdad \sigmae \k(\thetae)
  	\end{equation*}
  weakly in $L^2(0,T;\Hs)$, because, as we similarly did before, $\h$ and $\k$ are continuous and bounded by hypothesis~\ref{hyp:nonlinearities} and $\sigman$, $\thetan$ converges a.e. in $Q$ by \eqref{eq:limit_phin_ae}, \eqref{eq:limit_thetan_ae}. To conclude the proof, we only need to justify the fact that $\thetae$, $\phie$, and $\sigmae$ satisfy the initial conditions. On one hand, we observe that, along a subsequence, $\thetan \to \thetae$, $\phin \to \phie$, and $\sigman \to \sigmae$ strongly in $C^0([0,T];\Vs^*)$ at least. Thus, $\thetan(0) \to \thetae(0)$, $\phin(0) \to \phie(0)$, and $\sigman(0) \to \sigmae(0)$ strongly in $\Vs^*$. On the other hand, $\thetan(0) = P^n(\theta_0) \to \theta_0$, $\phin(0) = P^n(\phi_0) \to \phi_0$, and $\sigman(0) = P^n(\sigma_0) \to \sigma_0$ strongly in $\Hs$ by \eqref{eq:gal_limit_initial_datum_theta}--\eqref{eq:gal_limit_initial_datum_sigma}. By the uniqueness of the limit, we have
  	\begin{equation*}
  		\thetae(0) = \theta_0, \quad \phie(0) = \phi_0, \quad \sigmae(0) = \sigma_0.
  	\end{equation*}

\subsection{A priori estimate uniform in \(\varepsilon\)}
  	 Notice that the first and the second a priori estimates and their consequences at the Galerkin level are already independent of $\varepsilon$, so they pass to the limit as $n \to +\infty$ by lower semicontinuity, and are satisfied by $(\thetae,\phie, \mue, \sigmae)$. 

     The $\varepsilon$-dependence starts to appear in the third a priori estimate, and then, starting from it, it spreads out. Thus, we need to re-perform it, focusing on the problematic term, which now becomes
  	\begin{equation*}
  		\into\hbetae(\phi_{\varepsilon}(0)) \dx.
  	\end{equation*} 
  	Taking into account that $\betae$ can be bounded by $\hbeta$ (see the inequality \eqref{eq:yosida_property_hat}), and that $\hbeta(\phi_0)$ is integrable by assumption \eqref{hyp:initial_condition},  this integral is estimated  as follows
\begin{equation*}
	\into\hbetae(\phi_{\varepsilon}(0)) \dx	= \into\hbetae(\phi_0) \dx \leq \into\hbeta(\phi_0) \dx < +\infty,
\end{equation*}
independently from $\varepsilon$. Thus, also the third a priori estimate and its consequences hold with a constant independent of $\varepsilon$. By summarizing, we have that 
    \begin{equation}
        \begin{split}
            &\norm{\theta_{\varepsilon}}_{H^1(\Vs^*) \cap L^{\infty}(\Hs) \cap L^2(\Vs)} + \norm{\phi_{\varepsilon}}_{H^1(\Vs^*) \cap L^{\infty}(\Vs) \cap L^2(\Ws)} \\
            &\quad + \norm{\tau^{1/2}\phi_{\varepsilon}}_{H^1(\Hs)}+ \norm{\mu_{\varepsilon}}_{L^2(\Vs)} + \norm{\sigma_{\varepsilon}}_{H^1(\Hs) \cap L^{\infty}(\Vs) \cap L^2(\Ws)}\leq C
        \end{split}
        \label{pier-stima-eps}
    \end{equation}
for some constant $C$ depending on $M$ and on the problem's data, but independent of $\epsilon$.

\noindent Now we can derive an additional estimate, proceeding by comparison in the version of~\eqref{eq:ch2_def} written for $\betae$ (in place of $\beta$) and $(\thetae,\phie, \mue, \sigmae)$. In fact, in view of the regularity of the solution, we can write the equivalent equation
\begin{equation}
\tau \partial_t \phi_{\varepsilon} - \Delta \phi_{\varepsilon} + \betae(\phi_{\varepsilon}) + \pi(\phi_{\varepsilon}) - \chi \sigma_{\varepsilon} - \Lambda\theta_{\varepsilon}= \mu_{\varepsilon} 
\qquad \hbox{a.e. in } Q
\label{pier:ch2_def}
\end{equation}
and then comparing the terms in the light of \eqref{pier-stima-eps}. Hence, we can recover a uniform estimate for $\betae(\phi_{\varepsilon})$, namely
\begin{equation}
	\norm{\betae(\phie)}_{L^2(\Hs)} \leq C.
\label{pier-est-betaephie}
\end{equation}

\subsection{Passage to the limit as $\varepsilon \to 0$} From the estimates \eqref{pier-stima-eps}, \eqref{pier-est-betaephie} and standard compactness results, there exists a quadruplet $(\theta,\phi, \mu, \sigma)$ satisfying the regularity properties in~\Cref{defn:weak_sol} such that, as $\varepsilon \to 0$, along a nonrelabelled subsequence, the same convergences we had in \eqref{eq:limit_thetan_w}--\eqref{eq:limit_sigman_ae} hold, with the additional 
\begin{equation}\label{eq:limit_betae_w}
	\betae(\phie) \to \beta(\phi) \qquad \text{weakly} \qquad \text{in } L^2(0,T;\Hs).
\end{equation}
In fact, since $\betae(\phie)$ is uniformly bounded in $L^2(0,T;\Hs)$ it converges, along a subsequence, to a certain $\xi \in L^2(0,T;\Hs)$. Moreover, we know that $\phie \to \phi$ strongly in $L^2(0,T;\Hs)$. Thus, we may identify $\xi$ with $\beta(\phi)$ because of the strong-weak closeness of the graph of the maximal monotone operator $\beta$ (see \cite[][Proposition 2.5, p. 27]{brezis1973}). This is enough to pass to the limit in the approximate system, showing that the limit 
$(\theta,\phi, \mu, \sigma)$ is a weak solution to the PDE system~\eqref{eq:problem}--\eqref{eq:initial_datum}. Moreover, $(\theta,\phi, \mu, \sigma)$ satisfies the estimate~\eqref{est-pier1}.

\section{Regularity}\label{sect:reg}

This section contains the proof of~\Cref{thm:regularity}. To establish the existence of a more regular solution to our problem, we start again from the Faedo–Galerkin discretized system, taking advantage of the additional assumptions stated in \eqref{hyp:reg_initial_condition}. Then, we perform additional estimates and consequently obtain a limiting solution with the desired regularity.

\medskip

\noindent\textbf{Additional regularity of the discrete solution.} We aim to check that 
\begin{equation}\label{eq:phin_mun_add_reg}
	\phin \in H^2(0,T; \Vs^n), \quad \mun \in H^1(0,T; \Vs^n).
\end{equation}
We consider a fixed index $j = 0, \dots, n$. As already noticed, $\munj$ is a silent variable in the discrete Cahn--Hilliard equation \eqref{eq:ode_system_ch1}--\eqref{eq:ode_system_ch2}, and can be removed by substituting $\munj$, whose expression in given by equation \eqref{eq:ode_system_ch2}, into equation \eqref{eq:ode_system_ch1}. This way, we obtain:
\begin{equation}\label{eq:reg_phinj_derivative}
	\begin{split}
            (1+ \tau\gamma_j)\difft \phinj
            & = \intprod{(\lambdap \sigman - \lambdaa - \lambdae \thetan)\h(\phin)}{e_j} - \gamma_j^{\, 2}  \phinj \\
            &\quad-\gamma_j \intprod{\betae(\phin) + \pi(\phin)}{e_j} +\chi\gamma_j \sigmanj  + \Lambda \gamma_j\thetanj.
	\end{split}
\end{equation}
Looking at the right-hand side, we already know that $- \gamma_j^{\, 2}  \phinj+\chi\gamma_j \sigmanj  + \Lambda \gamma_j\thetanj$ belongs to the space $H^1(0,T)$. To conclude, since $e_j \in \Ws \subset L^{\infty}(\Omega)$, we only need to check that 
\begin{equation*}
	(\lambdap \sigman - \lambdaa - \lambdae \thetan)\h(\phin) -\gamma_j (\betae(\phin) + \pi(\phin)) \in H^1(0,T; L^1(\Omega)),
\end{equation*}
which is straightforward since $\h$, $\betae$, and $\pi$ are Lipschitz continuous, and $\sigman$, $\thetan$, and $\phin$ belong at least to $H^1(0,T;\Hs)$. Thus, $\phinj$ has the desired time-regularity and, consequently, $\phin$ satisfies \eqref{eq:phin_mun_add_reg}. By comparison in the equation \eqref{eq:ode_system_ch2}, it follows that $\munj \in H^1(0,T)$, so \eqref{eq:phin_mun_add_reg} is verified.\\

\noindent\textbf{Fourth priori estimate.} First of all, we make an observation.  In what follows, we will make use of the \emph{a priori} estimates derived in the proof of \Cref{thm:existence}. Recall that, at the Galerkin level, the first two estimates are uniform with respect to both $n$ and $\varepsilon$, whereas the third one is uniform only in $n$. This is due to the fact that we were unable to bound the term
\begin{equation*}
	\into \hbetae(\phi_{n,0}) \dx
\end{equation*}
uniformly in $\varepsilon$.
However, due to the stronger hypothesis on $\phi_0$, now we can do it. Since $\phi_0$ belongs to the space $\Ws\cap H^3(\Omega)$, it holds that (cf.~\eqref{pier-proiezione})
\begin{align}
    &\norm{\phi_{n,0}}_{\Ws} = \norm{P^n(\phi_0)}_{\Ws} \leq C \norm{\phi_0}_{\Ws}, \label{pier1}\\[2mm]
    & \norm{\Delta \phi_{n,0}}_{\Vs} = \norm{P^n(\Delta \phi_0)}_{\Vs}
    \leq C \norm{\Delta\phi_0}_{\Vs} \leq C \norm{\phi_0}_{\Ws\cap H^3(\Omega)}, \label {pier2}   
\end{align}
Thus, by the embedding $W\subset L^{\infty} (\Omega)$, it turns out that $\norm{\phi_{n,0}}_{L^{\infty}(\Omega)} \leq C$ for a constant $C$ that does not depend neither on $n$ nor $\varepsilon$. By the property in \eqref{eq:yosida_property_hat}, we infer that
\begin{equation*}
	\into \hbetae(\phi_{n,0}) \dx \leq \into \hbeta(\phi_{n,0}) \dx \leq C,
\end{equation*} 
since $\hbeta$ is continuous from hypothesis \ref{hyp:beta_pi}. This implies that, from now on, all the estimates previously 
derived are uniform with respect to both \( n \) and \( \varepsilon \). 
We point out another bound: in view of \eqref{pier-risolvente} and \eqref{pier-prop:betae}, it is not difficult to check that 
\begin{equation}
	\norm{\betae(\phi_{n,0})}_{\Vs} \leq C 
    \label{pier3}
\end{equation}
since $\, \nabla \betae(\phi_{n,0}) = \beta' (J_\varepsilon (\phi_{n,0}))
J_\varepsilon' (\phi_{n,0}) \nabla \phi_{n,0}\, $ a.e in $\Omega$, $\beta'$ is of class $C^1 $ and $J_\varepsilon$ is Lipschitz continuous with $J_\varepsilon (0)=0 $ and Lipschitz constant less than or equal to $1$.
Thanks to the bounds obtained in the third a priori estimate 
(see~\eqref{eq:gal_3rd_est}), we note that \eqref{eq:gal_dt_phi_mean_value_est} can be improved to
\begin{equation}
   	\norm{\overline{\partial_t\phin}}_{C^0([0,T])} \leq C \norm{(\lambdap \sigman - \lambdaa - \lambdae \thetan)\h(\phin)}_{L^\infty(\Hs)} \leq C.
    \label{pier8}
\end{equation}
Consequently, by testing equation~\eqref{eq:gal_ch2_def} with $v = |\Omega|^{-1}$ and making use of the estimates  \eqref{eq:gal_3rd_est} and \eqref{eq:gal_beta_est_infty}, we deduce that
\begin{equation}
\label{pier10}
    \begin{split}
        \norm{\overline{\mun}}_{C^0([0,T])} \leq &\tau \norm{\overline{\partial_t\phin}}_{C^0([0,T])}  + \sup_{t\in [0,T]} C\left( \norm{\betae(\phin(t))}_{L^1(\Omega)} + \norm{\phin(t)}_{L^1(\Omega)} + 1 \right)\\
        &+\sup_{t\in [0,T]} C\left(\norm{\sigman(t)}_{L^1(\Omega)} +  \norm{\thetan(t)}_{L^1(\Omega)} \right)  \leq C.
    \end{split}
\end{equation}

\noindent We now test equation~\eqref{eq:gal_temperature_def} with $ v = \partial_t \thetan $, equation~\eqref{eq:gal_ch1_def} with \( v = \partial_t \mun \), and the time-differentiated form of equation~\eqref{eq:gal_ch2_def} with \( v = \partial_t \phin \). Summing up the resulting equalities and noting a cancellation, we obtain
\begin{equation*}
	\begin{split}
            &\into |\partial_t \thetan|^2 \dx +\frac{1}{2} \difft \into |\nabla \thetan|^2 \dx + \frac{1}{2} \difft \into |\nabla \mun|^2 \dx\\
            &\quad \quad + \frac{1}{2} \difft \into |\tau^{1/2} \partial_t \phin|^2 \dx + \into |\nabla (\partial_t \phin)|^2 \dx + \into \betae'(\phin)|\partial_t \phin|^2 \dx\\
            &\quad= \into u \partial_t \thetan \dx 
            + \into (\ell - \Lambda)\partial_t \thetan  \partial_t \phin
            + \into (\lambdap \sigman - \lambdaa - \lambdae \thetan)\h(\phin) \partial_t \mun \dx\\
            &\quad \quad - \into \pi'(\phin)|\partial_t \phin|^2 \dx + \chi \into \partial_t \sigman \partial_t \phin \dx .
	\end{split}
\end{equation*}
Now, since $\betae$ is monotone and Lipschitz continuous, its derivative $\betae'$ is nonnegative, and the corresponding term on the left-hand side can be neglected. Then, using Young's inequality and the hypothesis~\ref{hyp:beta_pi} according to which $\pi'$ is bounded, we arrive at 
\begin{equation*}
	\begin{split}
            &\into |\partial_t \thetan|^2 \dx +\frac{1}{2} \difft \into |\nabla \thetan|^2 \dx + \frac{1}{2} \difft \into |\nabla \mun|^2 \dx\\
            &\quad \quad + \frac{1}{2} \difft \into |\tau^{1/2} 
            \partial_t \phin|^2 \dx + \into |\nabla (\partial_t \phin)|^2 \dx\\
            &\quad \leq C \into u^2 \dx + \frac{1}{2} \into |\partial_t \thetan|^2 \dx  + C \into |\partial_t \phin|^2 \dx\\
            &\quad \quad + C \into |\partial_t \sigman|^2 \dx + \into (\lambdap \sigman - \lambdaa - \lambdae \thetan)\h(\phin) \partial_t \mun \dx.
	\end{split}
\end{equation*}
Integrating in time over the interval $(0,t)$ for $t \leq T$ yields
\begin{equation}\label{eq:reg_a_priori_est_1}
	\begin{split}
            &\frac{1}{2}\bigg(\into |\nabla \thetan|^2 \dx + \into |\nabla \mun|^2 \dx + \into |\tau^{1/2} \partial_t \phin|^2 \dx \bigg)\\
            &\quad \quad + \frac{1}{2}\intto |\partial_t \thetan|^2 \dx \ds+ \intto |\nabla(\partial_t \phin)|^2 \dx \ds\\
            &\quad \leq \frac{1}{2}\left(\norm{\theta_{0,n}}_{\Vs}^2 + \norm{\nabla \mun (0)}_{\Hs}^2 + \norm{\tau^{1/2}\partial_t \phin(0)}_{\Hs}^2 \right) \\
            &\quad \quad + C\intto u^2 \dx\ds + C \intto |\partial_t \phin|^2 \dx \ds + C \into |\partial_t \sigman|^2 \dx \\
            &\quad \quad +  \intto (\lambdap \sigman - \lambdaa - \lambdae \thetan)\h(\phin) \partial_t \mun \dx \ds\\
            &\quad \eqcolon D_{0,n} + C \intto u^2 \dx \ds + I_1 + I_2 + I_3.
	\end{split}
\end{equation}
Concerning the terms on the last line, we observe that the constant $D_{0,n}$ accounts for the contribution from the initial data. Our goal is to show that $D_{0,n}$ is bounded independently of $n$. To this end, we recall that from \eqref{pier8} and \eqref{pier10} it follows that
\begin{equation}
   	|\overline{\partial_t\phin} (0) | + |\overline{\mun} (0) | 
   \leq C.
    \label{pier4}
\end{equation} 

\noindent Next, we consider equations~\eqref{eq:gal_ch1_def} and~\eqref{eq:gal_ch2_def} at the initial time. We subtract the mean value from~\eqref{eq:gal_ch1_def} and test the resulting equation with $\mathcal{N}(\partial_t\phin(0) - \overline{\partial_t\phin}(0))$. Then, we test~\eqref{eq:gal_ch2_def} with $(\partial_t\phin(0) - \overline{\partial_t\phin}(0))$ and add the two resulting expressions. By the property~\eqref{enum:eqN}, we observe that the terms involving $\mun(0)$ cancel out, as do all scalar products between mean values and either $\mathcal{N}(\partial_t\phin(0) - \overline{\partial_t\phin}(0))$ or $(\partial_t\phin(0) - \overline{\partial_t\phin}(0))$. Performing an integration by parts, we obtain
\begin{equation*}
    \begin{split}
        &\into |\nabla \N(\partial_t\phin (0) - \overline{\partial_t\phin} (0))|^2 \dx  + \tau \into |\partial_t\phin (0) - \overline{\partial_t\phin} (0)|^2 \dx\\
        &\quad = \dualprod{\lambdap \sigma_{0,n} - \lambdaa - \lambdae \theta_{0,n})\h(\phi_{0,n})}{\N(\partial_t\phin (0) - \overline{\partial_t\phin} (0))}\\
        & \quad \quad - \into (-\Delta \phi_{0,n} +      
        \betae ( \phi_{0,n} ) + \pi(\phi_{0,n}) - \chi \sigma_{0,n} -\Lambda \theta_{0,n}) (\partial_t\phin (0) - \overline{\partial_t\phin} (0))\dx .
    \end{split}
\end{equation*}     
The last term can be rewritten, again using~\eqref{enum:eqN}, as 
\begin{equation*}
- \into \nabla (-\Delta \phi_{0,n} +      
        \betae ( \phi_{0,n} ) + \pi(\phi_{0,n}) - \chi \sigma_{0,n} -\Lambda \theta_{0,n}) \cdot \nabla \N(\partial_t\phin (0) - \overline{\partial_t\phin} (0))\dx.
\end{equation*}
Then, by applying the Poincaré and Young inequalities, we conclude that
\begin{equation}
    \frac12 \into |\nabla \N(\partial_t\phin (0) - \overline{\partial_t\phin} (0))|^2 \dx  + \tau \into |\partial_t\phin (0) - \overline{\partial_t\phin} (0)|^2 \dx \leq C,
    \label{pier6}
\end{equation}
thanks to~\eqref{pier1}--\eqref{pier3} and the bounds on the initial data, which are under control due to the regularity assumptions $\sigma_0 \in \Vs$ and $\theta_0 \in \Vs$ (cf.~\eqref{hyp:reg_initial_condition} and~\eqref{pier-proiezione}). We point out that \eqref{pier6} and \eqref{pier4} yield in particular that $\norm{\partial_t\phin (0) }_{V^*}^2 \leq C$, due to the equivalence of norms. 

\noindent Finally, we take $v = \mun(0) - \overline{\mun}(0)$ in equation~\eqref{eq:gal_ch1_def} at the initial time. By carefully handling the terms and using the Poincaré and Young inequalities once more, we easily deduce that $\norm{\nabla \mun (0)}_{\Hs}^2 \leq C$, and hence, we ultimately obtain
\begin{equation}
D_{0,n}=\frac{1}{2}\left(\norm{\theta_{0,n}}_{\Vs}^2 + \norm{\nabla \mun (0)}_{\Hs}^2 + \norm{\tau^{1/2}\partial_t \phin(0)}_{\Hs}^2 \right) \leq C
\label{pier9}
\end{equation}
as desired.

\noindent The next step consists of estimating from above the last three terms on the right-hand side of \eqref{eq:reg_a_priori_est_1}. We deal with $I_1$ with the Ehrling's Lemma applied to the spaces $\Vs \subset \subset \Hs \subset \Vs^*$ and use the estimates \eqref{eq:gal_dt_phin_dualV} and \eqref{eq:gal_dt_phi_mean_value_est}, as follows:
\begin{equation}\label{pier7}
	\begin{split}
		&I_1 = C \intto |\partial_t \phin|^2 \dx \ds  \leq C \intt \norm{\partial_t \phin - \overline{\partial_t \phin}}_{\Hs}^2 \ds + C \intt |\overline{\partial_t \phin}|^2 \ds \\
		&\quad \leq \frac{1}{4} \intt \norm{\nabla (\partial_t \phin)}_{\Hs}^2 \ds + C \intt \norm{\partial_t \phin - \overline{\partial_t \phin}}_{\Vs^*}^2 \ds + C \intt |\overline{\partial_t \phin}|^2 \ds\\
		&\quad \leq  \frac{1}{4}  \intt \norm{\nabla(\partial_t \phin)}_{\Hs}^2 \ds + C.
	\end{split}
\end{equation} 
Moreover, we have that $$I_2= C \into |\partial_t \sigman|^2 \dx \leq C$$ 
by \eqref{eq:gal_3rd_est}.
We handle $I_3$ by integrating by parts
\begin{equation}
	\begin{split}
		I_3 =& -\intto (\lambdap \partial_t \sigman - \lambdae \partial_t \thetan)\h(\phin) \mun \dx \ds\\
		&- \intto (\lambdap \sigman - \lambdaa - \lambdae \thetan) \h'(\phin) \partial_t \phin \mun \dx \ds\\
		&+ \into (\lambdap \sigman - \lambdaa - \lambdae \thetan) \h(\phin) \mun \dx\\
		& - \into (\lambdap \sigma_{0,n} - \lambdaa - \lambdae \theta_{0,n}) \h(\phi_{0,n}) \mun (0) \dx \eqcolon I_{3,1}     + I_{3,2} + I_{3,3} + I_{3,4},
	\end{split}
\end{equation}
and then analyzing one by one its addends. The first one can be handled by Young's inequality, and the estimates \eqref{eq:gal_3rd_est}, \eqref{eq:gal_mu_est}:
\begin{equation}
	\begin{split}
		I_{3,1} &\leq C \intto (|\partial_t \sigman| + |\partial_t \thetan|) |\mun| \dx \ds\\
		&\leq \frac{1}{4} \intto |\partial_t\thetan|^2 \dx \ds + C \left(\intto |\partial_t \sigman|^2 \dx \ds  + \intto |\mun|^2 \dx \ds\right)\\
		& \leq \frac{1}{4} \intto |\partial_t\thetan|^2 \dx \ds + C. 
	\end{split}
\end{equation}
Regarding the second one, we employ the H\"older inequality, the estimates \eqref{eq:gal_3rd_est}, \eqref{eq:gal_mu_est}, and the Young inequality, obtaining: 
\begin{equation*}
	\begin{split}
		I_{3,2} &\leq C \intto (|\sigman| + |\thetan| + 1)|\partial_t \phin| |\mun| \dx \ds\\
		& \leq C \intt (\norm{\sigman}_{\Hs} + \norm{\thetan}_{\Hs} + 1) \norm{\partial_t \phin}_{L^4(\Omega)} \norm{\mun}_{L^4(\Omega)} \ds\\
		&\leq C \intt  \norm{\partial_t \phin}_{L^4(\Omega)}^2 \ds + C\intt \norm{\mun}_{\Vs}^2 \ds  \leq C \intt  \norm{\partial_t \phin}_{L^4(\Omega)}^2 \ds + C.
	\end{split}
\end{equation*}
To conclude, we add and subtract to $\partial_t \phin$ its mean value. Then, we apply Ehrling's lemma to the compact embeddings $\Vs \subset\subset L^4(\Omega) \subset \Vs^*$, similarly to what we did in~\eqref{pier7}. We deduce that
\begin{equation}\label{eq:reg_I12_ehrling_L4}
	\begin{split}
		I_{3,2} \leq C \intt  \norm{\partial_t \phin}_{L^4(\Omega)}^2 \ds + C  \leq  \frac{1}{4}  \intt \norm{\nabla(\partial_t \phin)}_{\Hs}^2 \ds + C.
	\end{split}
\end{equation}
We turn our attention to  $I_{3,3}$. Exploiting the H\"older and Poincar\'e inequalities, and the estimates~\eqref{eq:gal_3rd_est} and  \eqref{pier10}, we find
\begin{equation}
	\begin{split}
            I_{3,3} &\leq C \into (|\sigman| + |\thetan| + 1)|\mun| \dx \leq C \left( \norm{\sigman}_{\Hs} + \norm{\thetan}_{\Hs} + 1\right)\norm{\mun}_{\Hs}\\
            & \leq C \norm{\mun}_{\Hs} \leq C \norm{\nabla \mun}_{\Hs} + C|\omun| \leq \frac{1}{4} \norm{\nabla \mun}_{\Hs}^2 + C.	
        \end{split}
\end{equation}
Finally, recalling \eqref{pier9} and \eqref{pier4}, it is easy to conclude that $I_{3,4} \leq C$.

\noindent Now, going back to~\eqref{eq:reg_a_priori_est_1} and collecting all the intermediate estimates, upon rearranging the terms and adjusting the constants, we infer that
\begin{equation}
	\begin{split}
		&\into |\nabla \thetan|^2 \dx + \into |\nabla \mun|^2 \dx + \into |\tau^{1/2} \partial_t \phin|^2 \dx \\
		&\quad \quad + \intto |\partial_t \thetan|^2 \dx \ds + \intto |\nabla (\partial_t \phin)|^2 \dx \ds\leq  C,
	\end{split}
\end{equation}
whence
\begin{equation} \label{pier11}
	\begin{split}
		&\norm{\thetan}_{H^{1}(\Hs) \cap L^\infty(\Vs)} + \norm{\phin}_{H^{1}(\Vs) \cap L^{\infty}(\Vs)}+ \norm{\tau^{1/2}\phin}_{W^{1,\infty}(\Hs)} + \norm{\mun}_{L^\infty(\Vs)} \leq C.
	\end{split}
\end{equation}

\noindent \textbf{Fourth a priori estimate consequences.}   
     Taking $v = -\Delta \phin$ in \eqref{eq:gal_ch2_def}, which is admissible in our Faedo--Galerkin scheme, and exploiting integration by parts and monotonicity of $\betae$, it is straightforward to deduce that
 \begin{equation*}
\into | - \Delta \phin |^2 \dx  \leq \norm{ \mun - \tau \partial_t \phin - \pi(\phin) + \chi \sigman + \Lambda\thetan}_{\Hs}^2,
 \end{equation*}
with the right-hand side that is uniformly bounded in $L^\infty (0,T)$ due to~\eqref{pier11}. Then, from~\eqref{pier11} and elliptic regularity 
(see, e.g., \cite{Dautray_Lions_92}, \cite{Lions61}) it follows that 
\begin{equation} \label{pier12}
	\begin{split}
		\norm{\phin}_{L^{\infty}(\Ws)} \leq C.
	\end{split}
\end{equation} 
A similar procedure can be applied to \eqref{eq:gal_temperature_def} with the choice~$v = -\Delta \thetan$, in order to infer~that
\begin{equation} \label{pier13}
        \norm{\thetan}_{L^2(\Ws)} \leq C.
\end{equation} 
On the other hand, arguing as in~\eqref{pier14} we find out that
\[
\norm{\partial_t \phin}_{\Vs^*} \leq C \left(\norm{\nabla \mun}_{\Hs} + \norm{\sigman}_{\Hs} + \norm{\thetan}_{\Hs} + 1\right), 
\]
which leads to the estimate
\begin{equation} \label{pier15}
	\begin{split}
		\norm{\phin}_{W^{1,\infty}(\Vs^*)} \leq C.
	\end{split}
\end{equation} 
Moreover, if the choose $v = - \Delta \mun$ in equation \eqref{eq:gal_ch1_def}, we deduce that
\begin{equation*}
	\norm{-\Delta \mun}_{\Hs} \leq C \left(\norm{\partial_t \phin}_{\Hs} + \norm{\sigman}_{\Hs} + \norm{\thetan}_{\Hs} + 1\right) \leq C
\end{equation*}
and \eqref{pier11} and standard elliptic regularity results ensure that
\begin{equation} \label {pier16}
	\norm{\mun}_{L^2(\Ws)} \leq C.
\end{equation}

\noindent \textbf{Passages to the limit.} In the light of \eqref{pier11}--\eqref{pier16}, we can conclude that the limit $(\thetae, \phie, \mue, \sigmae)$ we found in the proof of the existence theorem as $n\to \infty$, enjoying the convergences \eqref{eq:limit_thetan_w}--\eqref{eq:limit_sigman_ae}, satisfies the additional estimate
\begin{equation} \label{pier17}
	\begin{split}
            &\norm{\thetae}_{H^{1}(\Hs) \cap L^\infty(\Vs)
            \cap L^2 (\Ws) } + \norm{\phie}_{W^{1,\infty}(\Vs^*)\cap H^{1}(\Vs) \cap L^{\infty}(\Ws)}\\
            &\quad  + \norm{\tau^{1/2}\phie}_{W^{1,\infty}(\Hs)} + \norm{\mue}_{L^\infty(\Vs)\cap L^2 (\Ws)} \leq C
	\end{split}
\end{equation}
which passes to the limit because it is independent of $\varepsilon$.
Note that, by the Sobolev embedding $W\subset L^\infty (\Omega)$,
\eqref{pier17} entails as well that $\norm{\phie}_{ L^{\infty}(Q)} \leq C$. Moreover, by comparison in the equation 
\eqref{pier:ch2_def} we find that
\begin{equation} \label{pier18}
	\norm{\betae(\phie)}_{L^{\infty}(\Hs)} \leq C.
\end{equation}
Passing to the limit as $\varepsilon \to 0$, all the estimates are preserved and therefore satisfied by the limit $(\theta, \phi, \mu, \sigma)$. Note that $\betae(\phie) \to \beta(\phi) $ weakly-$\ast$ in  $L^\infty(0,T;\Hs)$ due to \eqref{pier18}. 
Finally, we observe that, since the right-hand side in the equation
(cf.~\eqref{eq:temperature} and \eqref{hyp:control}) 
\[ \partial_t\theta  - \Delta \theta = u - \ell\partial_t \phi  \]
lies in $L^2(0,T; L^6(\Omega))$, and $\theta_0$ belongs to $L^{\infty}(\Omega)$, then it turns out that 
\begin{equation*}
	\norm{\theta}_{L^{\infty}(Q)} \leq C
\end{equation*}
by maximal parabolic regularity (see \cite[Chapter~III, Theorem~7.1, 
p.~181]{Ladyzenskaja_69}). This concludes the proof 
of~\Cref{thm:regularity}.

\section{Continuous dependence}\label{sect:cont_dep}
\noindent In order to prove \Cref{thm:continuous_dependence}, we consider two pairs $\{(\theta_i, \phi_i, \mu_i, \sigma_i)\}_{i=1,2}$ of strong solutions corresponding to  the initial data $\{(\theta_{0,i}, \phi_{0,i}, \sigma_{0,i})\}_{i=1,2}$ and to the assigned functions $\{u_i\}_{i=1,2}$. For convenience, in the following, we will employ the shorter notation
\begin{alignat*}{4}
	&\theta \coloneqq \theta_1 - \theta_2, \quad && \phi \coloneqq \phi_1 - \phi_2, && \mu \coloneqq \mu_1 - \mu_2, \quad && \sigma \coloneqq \sigma_{1} - \sigma_2,\\
	&\theta_0 \coloneqq \theta_{0,1} - \theta_{0,2}, \quad  && \phi_0 \coloneqq \phi_{0,1} - \phi_{0,2},  \quad && \sigma_0 \coloneqq \sigma_{0,1} - \sigma_{0,2}, \quad && u \coloneqq u_1 -u_2.
\end{alignat*}
Moreover, we recall the notation $f = \beta + \pi$ that we are going to use from now on.
First of all, we observe that  $(\theta, \phi, \mu, \sigma)$ satisfies  
\begin{subequations}\label{eq:difference_problem_eq_with_spaces}
	\begin{align}
            &\theta +\ell \phi - \Delta(1*_t\theta) =  \theta_{0}  +\ell \phi_{0} + (1*_t u), \label{eq:temperature_diff_def}\\
            &\partial_t \phi - \Delta \mu  =  (\lambdap \sigma - \lambdae \theta)\h(\phi_1) + (\lambdap \sigma_2 - \lambdaa - \lambdae theta_2) (\h(\phi_1) - \h(\phi_2)), \label{eq:ch1_diff_def}\\
            & \tau \partial_t \phi -\Delta \phi + f(\phi_1) -f(\phi_2) - \chi \sigma - \Lambda\theta  = \mu, \label{eq:ch2_diff_def}\\
            &\begin{aligned}
                \partial_t\sigma - \Delta(\sigma-\chi \phi) + \lambdab \sigma  = &  -\lambdac \sigma \h(\phi_1)- \lambdac \sigma_2 (\h(\phi_1) - \h(\phi_2))\\
                &- \lambdad \sigma \k(\theta_1) - \lambdad \sigma_2    (\k(\theta_1) - \k(\theta_2)),\label{eq:nutrient_diff_def}		
            \end{aligned} 
	\end{align}
\end{subequations}
a.e. in $Q$. This system is obtained straightforwardly by taking the difference of the systems (in the strong form) satisfied by $\{(\theta_i, \phi_i, \mu_i, \sigma_i)\}_{i=1,2}$, and integrating the temperature equation in time. We multiply \eqref{eq:temperature_diff_def} by $ R \theta$, where $R$ is a positive (big) constant yet to be determined, \eqref{eq:ch1_diff_def} by $\phi$, \eqref{eq:ch2_diff_def} by $-\Delta \phi$, and \eqref{eq:nutrient_diff_def} by $\sigma$. We sum all these equalities and integrate over $\Omega$, finding 
\begin{equation}\label{eq:cnt_dep_ineq_1}
	\begin{split}
		& R\into |\theta|^2 \dx + \frac{R}{2} \difft \into |\nabla(1 *_t \theta)|^2 \dx + \frac{1}{2} \difft \into |\phi|^2 \dx + \frac{1}{2} \difft \into |\tau^{1/2} \nabla \phi|^2 \dx\\
		&\quad \quad + \into |-\Delta \phi|^2 \dx + \frac{1}{2} \difft \into |\sigma|^2 \dx + \into |\nabla \sigma|^2 \dx + \lambdab \into |\sigma|^2 \dx\\
		& \quad =  R\into (\theta_{0} + \ell \phi_{0})\theta \dx + R \into (1*_t u)\theta \dx \\
		&\quad \quad- R \ell \into \phi \theta \dx + \into (\lambdap \sigma - \lambdae \theta)\h(\phi_1)\phi \dx\\
		&\quad \quad  + \into  (\lambdap \sigma_2 - \lambdaa - \lambdae \theta_2) (\h(\phi_1) - \h(\phi_2)) \phi \dx\\
		&\quad \quad  - \into (f(\phi_1) - f(\phi_2)) (-\Delta\phi) \dx + 2 \chi \into \sigma (-\Delta \phi) \dx +\Lambda\into \theta (-\Delta \phi) \dx\\
		&\quad \quad + \into \left[ -\lambdac \sigma \h(\phi_1)- \lambdac \sigma_2 (\h(\phi_1) - \h(\phi_2))\right] \sigma \dx\\
		&\quad \quad + \into \left[- \lambdad \sigma \k(\theta_1) - \lambdad \sigma_2    (\k(\theta_1) - \k(\theta_2)) \right] \sigma \dx.
	\end{split}
\end{equation}
Then, recalling assumptions~\ref{hyp:nonlinearities}, \ref{hyp:beta_pi}, and the regularity estimate~\eqref{pier21}, we exploit the following facts:
\begin{itemize}
    \item $\h$  and $\k$ are bounded and Lipschitz continuous;
    \item $\theta_2$ is bounded in $L^{\infty}(Q)$;
    \item $\sigma_2$ is bounded in $L^2(0,T; \Ws)$, and it holds that $\|\sigma_2\|_{L^\infty(\Omega)} \leq C\|\sigma_2\|_{W}$ a.e.\ in $(0,T)$;
    \item $f$ is locally Lipschitz continuous, and $\phi_1$, $\phi_2$ are bounded in $L^{\infty}(Q)$.
\end{itemize}
Using Young's inequality and the bounds above, the right-hand side of~\eqref{eq:cnt_dep_ineq_1} can be estimated by
 \begin{equation}\label{eq:cnt_dep_ineq_2}
 	\begin{split}
            &\frac{R}2 \into |\theta|^2 \dx  + C\/R \left( \norm{\theta_0}_{\Hs}^2 +\norm{\phi_0}_{\Hs}^2 +\into |1*_t u|^2 \dx + \into |\phi|^2 \dx \right)\\
            &\quad+ \frac12\into |\theta|^2 \dx +  C\left( \into |\sigma|^2 \dx + \into |\phi|^2 \dx \right) + C\left(1+\norm{\sigma_2}_{\Ws}\right)\into |\phi|^2 \dx\\
            &\quad + 8 \chi^2\into |\sigma|^2 \dx +2 \Lambda^2 \into |\theta|^2 \dx + \frac{1}{2} \into |-\Delta \phi|^2 \dx \\
            &\quad+ C \into |\sigma|^2 \dx + C \norm{\sigma_2}_{\Ws} \into |\phi|^2+ |\sigma|^2 \dx  + C          \norm{\sigma_2}_{\Ws}\into |\theta|\/ |\sigma| \dx.
 	\end{split}
 \end{equation}
Let us briefly comment on the fact that, at this stage, the constants $C$ used above depend on the norms of the solutions \( \{(\theta_i, \phi_i, \mu_i, \sigma_i)\}_{i=1,2} \) through the estimate~\eqref{pier21}. However, as a consequence of the proof we are currently carrying out, we will establish the uniqueness of the solution. This implies that any solution must coincide with the one constructed in \Cref{thm:existence,thm:regularity}. Therefore, it satisfies the a priori estimates derived therein, which ensure that the norm of each component of the solution can be bounded by a constant depending only on the data of the problem---such as the final time $T$, the domain $\Omega$, the constant $R$, the initial data... As a result, the same conclusion applies to the constants $C$ appearing above.
Most of the terms in the previous expression \eqref{eq:cnt_dep_ineq_2} do not need any further treatment because, after time integration, we are going to apply the Gronwall Lemma. Only the last one requires some additional calculations. By the Young inequality, we find
\begin{equation}\label{eq:cnt_dep_ineq_3}
	\begin{split}
		C \norm{\sigma_2}_{\Ws}\into |\theta|\/ |\sigma| \dx
		\leq \frac{1}{2} \into |\theta|^2 \dx + C\norm{\sigma_2}_{\Ws}^2 \into |\sigma|^2 \dx.
	\end{split}
\end{equation}
Now we collect the contibutions \eqref{eq:cnt_dep_ineq_1}--\eqref{eq:cnt_dep_ineq_3}, fix $R > \frac R 2 + 1 + 2\Lambda^2$ (cf.~the coefficients of the terms with $\theta $ on the right-hand side), then move to the left-hand side the terms in $\theta$ and $-\Delta\phi$, and integrate in time over $(0,t)$, obtaining:
\begin{equation*}
	\begin{split}
		& \into |\nabla(1 *_t \theta)|^2 \dx +  \into |\phi|^2 \dx + \into |\tau^{1/2} \nabla \phi|^2 \dx + \into |\sigma|^2 \dx\\
		&\quad \quad + \intto |\theta|^2 \dx\ds  + \intto |-\Delta \phi|^2 \dx\ds  + \intto |\sigman|^2\dx \ds + \intto |\nabla \sigman|^2\dx \ds\\
		& \quad \leq C \bigg[\norm{\theta_0}_{\Hs}^2  + \norm{\phi_0}_{\Hs}^2  +\tau^{1/2} \norm{\nabla\phi_0}_{\Hs}^2  +\norm{\sigma_0}_{\Hs}^2 +  \intto |1*_t u|^2 \dx \ds \\
		&\quad \quad+  C \intt \left(1+\norm{\sigma_2}_{\Ws}\right)\into |\phi|^2 \dx \ds  + C\intt\bigl(1+\norm{\sigma_2}_{\Ws}^2\bigr) \into |\sigma|^2 \dx \ds \bigg].
	\end{split}
\end{equation*}
Hence, since the function $t\mapsto 1+\norm{\sigma_2}_{\Ws}^2 $ is bounded in $L^1(0,T)$, by the Gronwall Lemma we infer that
\begin{equation}\label{eq:cnt_dep_est_1}
    \begin{split}
        &\norm{\theta}_{L^2(\Hs)}^2 + \norm{\nabla(1*_t\theta)}_{L^{\infty}(\Hs)}^2 + \norm{\phi}_{L^{\infty}(\Hs)            \cap L^2(\Ws)}^2\\
        &\quad \quad + \norm{\tau^{1/2} \phi}_{L^{\infty}(\Vs)}^2 + \norm{\sigma}_{L^{\infty}(\Hs) \cap L^2(\Vs)}^2\\
        &\quad \leq C\left( \norm{\theta_0}_{\Hs}^2 + \norm{\phi_0}_{\Hs}^2 +  \norm{\tau^{1/2}\phi_0}_{\Vs}^2  +  \norm{\sigma_0}_{\Hs}^2 +  \norm{1*_t u}_{L^2(\Hs)}^2 \right). 
    \end{split}
\end{equation}
Finally, we want to improve this result for the temperature variable $\theta$. To do so, we consider the equation satisfied by $\theta$ a.e. in $Q$, which is the following
(cf. \eqref{eq:temperature})
\begin{equation*}
    \partial_t (\theta +\ell\phi) - \Delta \theta = u.
\end{equation*}
We multiply it by $\theta +\ell\phi$ and integrate over $\Omega \times (0,t)$, finding
\begin{equation*}
    \begin{split}
        &\frac{1}{2}\into |\theta +\ell\phi|^2 \dx + \intto |\nabla \theta|^2 \dx \ds\\
        &\quad = \frac{1}{2} \into |\theta_0 +\ell\phi_0|^2 \dx + \intto u (\theta +\ell\phi) \dx \ds  - \ell\intto \nabla \theta \cdot \nabla \phi \dx \ds. 
    \end{split}
\end{equation*}
We apply the Young inequality and the previously proved inequality \eqref{eq:cnt_dep_est_1} to treat the terms on the right-hand side, leading to
\begin{equation*}
    \begin{split}
        &\frac{1}{2} \into |\theta +\ell\phi|^2 \dx + \frac{1}{2} \intto |\nabla \theta|^2 \dx \ds \\
        &\quad \leq C \bigg(\norm{\theta_0}_{\Hs}^2  + \norm{\phi_0}_{\Hs}^2+ 
        \intto |u|^2 \dx \ds \\
        &\quad \quad +  \intto |\theta|^2 +|\phi|^2 \dx \ds + \intto |\nabla \phi|^2 \dx \ds \bigg)\\
        &\quad \leq C \bigg(\norm{\theta_0}_{\Hs}^2  + \norm{\phi_0}_{\Hs}^2  +\tau^{1/2} \norm{\nabla\phi_0}_{\Hs}^2  +\norm{\sigma_0}_{\Hs}^2  + \intTo |u|^2 \dx \ds \bigg).
    \end{split}
\end{equation*}
Then, we end up with
\begin{equation}\label{eq:cnt_dep_est_2}
    \begin{split}
        &\norm{\theta +\ell\phi}_{L^{\infty}(\Hs)} + \norm{\nabla \theta}_{L^2(\Hs)}\\
        &\quad \leq C\left(\norm{\theta_0}_{\Hs}^2 + \norm{\phi_0}_{\Hs}^2 +  \norm{\tau^{1/2}\phi_0}_{\Vs}^2  +  \norm{\sigma_0}_{\Hs}^2  +  \norm{u}_{L^2(\Hs)}^2 \right)
    \end{split}
\end{equation}
and the continuous dependence inequality~\eqref{pier22} follows easily from the estimates \eqref{eq:cnt_dep_est_1}--\eqref{eq:cnt_dep_est_2}. 
Therefore, \Cref{thm:continuous_dependence} is completely proved.

\section*{Acknowledgments}
The authors are members of GNAMPA (Gruppo Nazionale per l’Analisi Matematica, la Probabilità e le loro Applicazioni) of INdAM (Istituto Nazionale di Alta Matematica).~P.C. and E.R. gratefully mention some support from the Next Generation EU Project 
No.~P2022Z7ZAJ (A unitary mathematical framework for modelling muscular dystrophies) and the MIUR-PRIN Grant 2020F3NCPX “Mathematics for industry 4.0 (Math4I4)".
\printbibliography

\end{document}